\documentclass[final, authoryear, 3p,times, twocolumn]{elsarticle}

\usepackage{graphics}

\usepackage{amssymb}
\usepackage{amsmath}

\usepackage{lineno}

\begin{document}

\begin{frontmatter}



\title{Estimating numerical errors due to operator splitting in global atmospheric chemistry models: transport and chemistry.}


\author[label1,label2]{Mauricio Santillana\corref{cor1}}
\ead{msantill@fas.harvard.edu}
%
%

\author[label3]{Lin Zhang}
\ead{zhanglg@pku.edu.cn}
%
\author[label1]{Robert Yantosca}
\ead{yantosca@seas.harvard.edu}
%
\address[label1]{School of Engineering and Applied Sciences, Harvard University, Cambridge, MA 02138, United States}
\address[label2]{Computational Health Informatics Program, Boston Children's Hospital, Boston, MA 02115, United States}
\address[label3]{Laboratory for Climate and Ocean-Atmosphere Sciences, Department of Atmospheric and Oceanic Sciences, School of Physics, Peking University, Beijing, 100871, China}
\cortext[cor1]{Corresponding author. Tel.: 617 495-2891 }



\begin{abstract}
We present upper bounds for the numerical errors introduced when using operator splitting methods to integrate transport and non--linear chemistry processes in global chemical transport models (CTM). We show that (a) operator splitting strategies that evaluate the stiff non--linear chemistry operator at the end of the time step are more accurate, and (b) the results of numerical simulations that use different operator splitting strategies differ by at most 10\%, in a prototype one--dimensional non--linear chemistry-transport model. We find similar upper bounds in operator splitting numerical errors in global CTM simulations.


\end{abstract}


\end{frontmatter}


\section{Introduction}
\label{sec:intro}
Global tropospheric chemistry transport models (CTM) are used to address important issues ranging from air quality to climate change. In order to continuously improve their performance, it is of crucial importance to understand and quantify the diverse sources of uncertainties and errors present in them. We group these in three different categories,  (\textit{i}) errors and uncertainties coming from observations and data used in our models (such as emission inventories, wind fields, reaction rates); (\textit{ii})  errors coming from our choice of governing equations (or mathematical model), parametrizations, and the level of complexity of the physical modules included in our formulation; and (\textit{iii}) numerical errors coming from the choice of algorithms we use to solve the governing equations using computers \citep{bib:Ent02, bib:Zha11}. \\

In this study, we focus our attention on estimating the magnitude of numerical errors (\textit{iii}), in particular, those arising from the choice of operator splitting technique utilized 
to integrate in time the transport and chemistry operators in real-life global CTMs. In order to achieve this, we numerically extend the results introduced for the linear diffusion-reaction case in \citep{bib:Spo00}, to a non-linear 1-D chemistry-transport numerical model. The latter numerical results provide us with a framework to estimate upper bounds for operator splitting errors in the fully non-linear 3-D state-of-the-art global CTM: GEOS-Chem \citep{bib:Bey01}. To the best of our knowledge, our contribution is the first in estimating operator splitting errors in the context of real-life global atmospheric chemistry simulations. \\

CTMs simulate the dynamics of chemical species in the atmosphere by numerically integrating a set of coupled nonlinear partial differential equations  
of the type:

\begin{equation} 
\dfrac{\partial C_i}{\partial t}+\nabla\cdot \left({\boldsymbol u} \:C_i
\right)=\nabla\cdot \left( \rho K \nabla \dfrac{C_i}{\rho} \right)+ P_i(C_j)- C_i L_i(C_j)+Q_i - S_i 
\label{eq:advec-reac}
\end{equation}
for $i=1,...,N$; where $C_i({\boldsymbol x},t)$  represents the spatio-temporal evolution of the concentration of species $i$ (typically over a hundred species are considered), ${\boldsymbol u}({\boldsymbol x},t)$ is the wind velocity, 
$\rho$ is the air density, $K$ the eddy diffusivity matrix, $P_i$ are the nonlinear production terms, $L_i$ are the destruction terms, $Q_i$ are the volume emission sources, and $S_i$ are the sinks (ex. precipitation or in-cloud removal). See \cite{bib:Spo07} for a detailed description of these equations. \\

Due to the dimensions of grid boxes in global CTMs, like GEOS-Chem (with hundreds of kilometers in the horizontal versus tens to hundreds of meters in the vertical), intertial vertical transport processes in this global models are simulated (a) using vertical mass fluxes schemes that ensure that the horizontal air flow is divergent-free ($\nabla_{hor}\cdot \boldsymbol u=0$), (b) using convection parametrizations, and (c) using a boundary layer mixing algorithm \citep{bib:Lin96,bib:All96, bib:Wil06, bib:Pra08}. In addition, horizontal diffusion due to numerical errors in transport schemes are typically higher than their Eddy difusivity counterpart, as measured by aircraft missions \citep{bib:Pis09, bib:Wil06, bib:Ras07, bib:San13}. As a consequence, the first term of the right-hand side of equation (\ref{eq:advec-reac}), which models the dynamics of intertial vertical transport as an eddy diffusion process, is not explicitly integrated in global CTMs; and the governing equations (\ref{eq:advec-reac}) are sometimes written \citep{bib:Ras07, bib:San10, bib:San13} in a simplified way as
\begin{equation} 
\dfrac{\partial C_i}{\partial t}+{\boldsymbol u}\cdot \nabla C_i=P_i(C_j)- C_i L_i(C_j)+Q_i - S_i. 
\label{eq:advec-reac2}
\end{equation}

The chemistry operator on the right-hand-side of equations (\ref{eq:advec-reac2}) models the chemical interaction of atmospheric species whose lifetimes range from milliseconds to many years. The chemistry operator is very stiff as a consequence of this large range of time-scales  and thus, implicit-in-time methods are an appropriate choice to integrate equations (\ref{eq:advec-reac}). Traditional methods, such as the method of lines, aimed at achieving this task in realistic 3D simulations, involve solving for an enormous number of degrees of freedom at each time step in a coupled fashion ( $10^8 \approx$ 100 chemical species in $\sim 10^6$ grid cells, for a $1^{\circ}\times 1^{\circ}$ spatial resolution). This is due to the inter-species coupling in the chemistry operator and the spatial coupling in the transport operator. In practical situations, however, efficient computational algorithms to integrate equations (\ref{eq:advec-reac}) use operator splitting strategies that allow the explicit time--integration of the transport and implicit time--integration of the chemistry operators separately and sequentially, thus, reducing significantly the degrees of freedom solved in a coupled fashion at a given time step. This is done at the expense of a loss of accuracy in the approximate solution \citep{bib:HunVer03}. \\

Estimating the magnitude of the numerical errors introduced by the time--integration of equations (\ref{eq:advec-reac}) in realistic 3-D computer simulations is a hard task since no relevant analytic solution can be used as a reference to estimate them. In theory, estimates of  these errors depend directly on the regularity properties of the analytic solution of equations (\ref{eq:advec-reac}), the set of initial and boundary conditions, and the chosen numerical scheme \citep{bib:Guo86, bib:Ise96, bib:Ern04, bib:Bre08}. In this study, we assume that the analytic solution of equations (\ref{eq:advec-reac}) is unique and regular enough so that numerical error estimates can be expressed as inequalities of the form (\ref{eq:error_est}). 
Operator splitting errors, as well as numerical errors arising from the time--integration of the chemistry operator depend explicitly on the magnitude of the chosen time steps, while numerical errors coming from the time--integration of the transport operator depend both on the time step and on the grid size. This fact, in combination with an expression of the analytic solution of equations (\ref{eq:advec-reac}), is exploited to obtain the exact magnitude of operator splitting errors in our one-dimensional proto-type transport-chemistry numerical model. \\

Our one-dimensional numerical experiments show three main results: (a) operator splitting sequences where the stiff non--linear chemistry operator is evaluated at the end of the time step are more accurate than those where the transport is evaluated lastly, independently of the operator splitting time-step, as in the linear case introduced in \citep{bib:Spo00}; (b) the results of numerical simulations that use different operator splitting strategies differ by at most 10\%; and (c) numerical errors coming from the integration of the transport operator are much bigger than those coming form the operator splitting technique for spatial and temporal scales comparable to those used in global CTM. We use this fact, and evidence from papers such as \citep{bib:Wil06, bib:Ras07, bib:Pra08, bib:San13}, to suggest that in realistic 3D simulations, errors due to operator splitting are much smaller than those introduced by transport schemes. 

\section{Numerical error estimation}

Upper bounds of the numerical errors introduced by solving partial differential equations with regular boundary and initial conditions, using a given numerical scheme, can be expressed by inequalities represented as
\begin{equation}
||C(x,t)-C_h(x,t)||_{_{V_1}}\leq M_1\:\Delta t\:^\alpha  + M_2\: \Delta x\:^\beta 
\label{eq:error_est}
\end{equation}
where $C(x,t)$ is the true solution of the partial differential equation, $C_h(x,t)$ the numerical approximation, $\Delta t$ and $\Delta x$ are the time step and grid size respectively, $\alpha$ and $\beta$ are exponents (typically larger than one) that determine the order of convergence of the method in time and space respectively, $M_1$ and $M_2$ are constants that depend on the regularity of the true solution $C(x,t)$ and parameters in the equation, and $||\cdot||_{V_1}$ is the norm in the appropriate Banach space $V_1$. For a convergent method, as $\Delta t\rightarrow 0$ and $\Delta x\rightarrow 0$, the numerical error vanishes, ({\it i.e.} $||C-C_h||_{_{V_1}}\rightarrow 0$) and the numerical approximation $C_h$ converges to the true solution $C$, in the normed space $V_1$. More details about the integral representation (equation \ref{eq:error_est}) of numerical errors due to discretization of partial differential equations can be found in: \cite{bib:Guo86, bib:Ise96, bib:Ern04, bib:Bre08}\\

For the specific set of partial differential equations (\ref{eq:advec-reac}), operator splitting errors and errors coming from the numerical integration of the chemistry operator (where no coupling in space exists) contribute to the first term on the right-hand-side of inequality (\ref{eq:error_est}), whereas, numerical errors from the integration of the transport operator contribute to the first and second terms of the right-hand-side of inequality (\ref{eq:error_est}). Quantifying the independent contribution of each processes to each term of inequality (\ref{eq:error_est}) is not simple in practical applications. In the following section, we show how to estimate the magnitude of operator splitting errors in the absence of other numerical error coming from the time--integration of the transport and chemistry operators.   \\


\subsection{Operator splitting techniques and error estimation}
\label{sec:op_split}

Classical approaches to estimate the numerical errors introduced by operator splitting approaches are based on asymptotic expansions of exponential operators (linear case) and Lie operator formalism (nonlinear case). For completeness, we briefly describe important results of the linear analysis of operator splitting methods in this section. We refer the reader to \citet{bib:LanVer99, bib:Spo00, bib:HunVer03} and the references therein for more details. 
In this section, it is assumed that the time--integration of each operator separately can be found exactly giving rise to no numerical error, {\it i.e.} the numerical errors discussed below come only from the choice of the operator splitting technique. \\

We use as an example the linear evolution equation, 
\begin{equation}
\dfrac{d { v}}{d t}=Av+Bv, \qquad v(0)=v_0, \qquad v\in \mathbb{R}^n
\label{eq:op_split1}
\end{equation}
where $A$ and $B$ are linear operators. One of these operators could represent the linear spatial differential operator $d/dx$ (transport) in equations (\ref{eq:advec-reac}). The analytic solution for this problem is given by:
\begin{equation}v=v_0 \exp((A+B)t)
\label{eq:exact}
\end{equation} The simplest operator splitting method, called Godunov and denoted by $(A-B)$, can be obtained for $t\in[0,\Delta t]$ by solving the two evolution equations in sequence as:
\begin{equation}
\left\{
\begin{aligned}
\dfrac{d {v^*}}{d t} &=Av^*, \qquad  &v^*(0) = & v_0 \qquad & in \; [0, \Delta t]\\
\dfrac{d {v^{**}}}{d t}& =Bv^{**}, \qquad & v^{**}(0) = &v^*(\Delta t) \qquad & in \; [0, \Delta t].
\end{aligned}
\right.
\end{equation}
The value for $v$ at $t=\Delta t$ is given by $v_{AB}(\Delta t)=v^{**}(\Delta t)$. The solution obtained with this operator splitting method at $t=\Delta t$ is given by
\begin{equation}v_{AB}(\Delta t)=v_0 \exp(B \Delta t)\exp(A\Delta t)
\label{eq:godunov}
\end{equation} The exact solution (\ref{eq:exact}) and the solution $v_{AB}$ in the previous equation will be the same if 
\[\exp((A+B)\Delta t )=\exp(B \Delta t)\exp(A\Delta t).
\] This will happen if the operators $A$ and $B$ commute (think of matrices), \textit{i.e.} if $AB=BA$. When $AB\neq BA$, then the (point-wise) local-in-time numerical error associated to solving problem (\ref{eq:op_split1}) using Godunov's operator splitting technique can be shown to be 
\begin{equation}
le_{AB}=\dfrac{(AB-BA)}{2} \Delta t^2 v_0 
\label{eq:ab_error}
\end{equation}
which leads to a global error $\mathcal{O}(\Delta t)$, \textit{i.e.} $||v-v_{AB}||\leq M_{AB}\;\Delta t$ (for a constant $M_{AB}$ that depends only on the regularity of the analytic solution $v$). Since the numerical error vanishes as $\Delta t \rightarrow 0$, Godunov's method is a convergent first order method in time, in the linear case.  
Another simple Godunov operator splitting can be obtained by reversing the order of evaluation of the operators $A$ and $B$ to obtain the $(B-A)$ method ($v_{BA}$). A more accurate and symmetric operator splitting method, often referred to as Strang method (Strang, 1968), can be obtained by averaging the output of the two previous methods, i.e. $v_{S}(\Delta t)=\frac{1}{2}(v_{AB}+v_{BA})$. It can be shown that the Strang operator method is globally second order accurate, \textit{i.e.} $||v-v_{S}||\leq M_{S}\;\Delta t^2$ for a constant $M_S$ \citep{bib:Spo00, bib:HunVer03}.\\

The linear analysis presented above may fail and lead to different convergence results if one of the operators is stiff, {\it i.e.} if the dynamics of one operator take place in much faster time scales than the dynamics in the other operator \citep{bib:Spo00}.  This can be seen by introducing a small parameter $\epsilon$ (representing the ratio between fast time scales in the stiff operator and the slow time scales of the other operator) and re-writing the linear evolution equation (\ref{eq:op_split1}) as a singular perturbation equation by re-defining
\begin{equation}
A=\dfrac{\chi(\epsilon)}{\epsilon} \qquad \qquad \text{and} \qquad \qquad B=T.
\label{eq:defn_stiff}
\end{equation}
For our purposes, one can identify the chemistry operator with the stiff operator $\chi/\epsilon$, (the nonlinear chemistry can be, locally-in-time and space, approximated by a linear and stiff mechanism at least for some subset of fast species), and identify the transport operator with the slow operator $T$, for which the dynamics takes place in a more confined range of time scales (as represented by our global models). It is shown in \citep{bib:Spo00} that the local error for the $\left(\frac{\chi}{\epsilon}-T\right)$ Godunov method becomes (compare to equation (\ref{eq:ab_error})):
\begin{equation}
le_{\epsilon}\sim\dfrac{(\chi\: T - T\:\chi)}{\epsilon} \Delta t^2 v_0 \\
\label{eq:stiff_error}
\end{equation}
leading to a global error $\mathcal{O}\left(\frac{\Delta t}{\epsilon}\right)$, implying that $||v-v_{\epsilon}||\leq M_{\epsilon}\left(\frac{\Delta t}{\epsilon}\right)$. Note that convergence of the operator splitting method, in this case, can only be guaranteed provided the operator splitting time step, $\Delta t$, is small enough to satisfy $\Delta t \ll \epsilon$ so that higher order terms, $\mathcal{O}(\frac{\Delta t}{\epsilon})^k$, will indeed vanish as $k\rightarrow \infty$ in the Taylor expansion of the error. \\

In atmospheric chemistry simulations, we use operator splitting methods to integrate in-time two operators in equations  (\ref{eq:advec-reac}): transport and chemistry. Transport and chemistry are known to commute when the velocity field is divergent free and chemistry is independent of the spatial location. In real atmospheric situations, these conditions are typically not met. Indeed, the non-linear chemistry operator depends dynamically on the geographic location (due to photolysis), and atmospheric wind fields are in general not divergent-free. 

The result of the linear analysis above suggests that 
%
%
%
operator splitting approaches will converge only if the operator splitting time step is much smaller than the lifetime of the fastest species in the chemistry mechanism ($\Delta t \ll \epsilon$). This is also the criterion established to ensure stability and convergence of explicit-in-time chemistry solvers, and suggests the use of prohibitively small operator splitting time-steps in order to guarantee convergence of the method. In practice, however, the use of implicit schemes to integrate the chemistry operator in global chemistry models leads to the choice of large operator splitting time-steps compared to the intrinsic stiffness of the chemistry system ($\Delta t\gg \epsilon$). As a consequence, and according to expression (\ref{eq:stiff_error}), we may expect to observe large operator splitting errors when solving equations (\ref{eq:advec-reac}) with stiff and potentially non--linear chemistry operators.\\


It is argued in \citet{bib:Spo00}, that operator splitting errors, even in the presence of large operator splitting time steps (such that $\Delta t\gg \epsilon$), may not be as big as suggested by expression (\ref{eq:stiff_error}). \citet{bib:Spo00} argues that the stiffness of the system can be balanced by the existence of an underlying reduced model (low-dimensional manifold) describing the dynamics of the system and thus, by choosing the appropriate order of operator evaluation in a time-step, the splitting error may be bounded even with the increase of stiffness. Moreover, he shows for the linear case that sequences where the stiff operator is evaluated at the end of the time step lead to convergent and accurate methods in a one dimensional diffusion-chemistry toy example, even for large operator splitting time steps. In solving equations (\ref{eq:advec-reac}), examples of these sequences include: Transport--Chemistry and Chemistry--Transport--Chemistry. \\

Intuitively speaking, evaluating the transport operator at the end of the time step sets the state of the system far from the underlying low dimensional manifold driving the chemical system and provides an initial condition $v_0$ for the next time evaluation that enhances error propagation. This is avoided by evaluating the stiff chemistry operator at the end of the time step. The existence of these reduced models driving the dynamics in regional and global atmospheric chemistry models has been found in \citet{bib:LowTom00, bib:San10, bib:Ras07}, suggesting that the operator splitting order should be selected carefully. To the best of our knowledge, a careful investigation of these errors in the realistic non--linear case does not exist so far and thus we aim at achieving this here. \\


Isolating operator splitting errors in practical global atmospheric chemistry models is not straightforward, first, because we lack expressions for the analytic solution of the system in realistic circumstances, and second, since the solutions of the chemistry and transport operators, separately, are obtained using numerical schemes and thus are not exact as it was assumed in the previous analysis. In order to estimate upper bound estimates of operator splitting errors we proceeded as follows. We first found sharp estimates of numerical errors in a 1D non--linear chemistry-transport prototype problem with a known analytic solution. We designed this 1D problem to resemble the interaction of numerical errors in the time--integration of the transport and (stiff) non--linear chemistry, when using operator splitting methods, at spatial and time scales used in 3D global simulations. Our 1D findings guide our methodology to understand the differences observed between the outputs of 3D global simulations using different operator splitting strategies. We performed multiple 3D global simulations in order to further understand additional numerical errors, due to the time integration of relevant processes (emissions, convective transport, and deposition) inherently solved with operator splitting approaches.


\section{One-dimensional advection-reaction system}
\label{sec:one_dim_numerics}

We considered a one-dimensional advection-reaction system that can be solved analytically and thus exact values of numerical errors can be obtained. The system is characterized by a constant wind field throughout the domain, and a three-species ($NO$, $NO_2$ , $O_3$) stiff non--linear chemistry-mechanism modeling the $NO_x(NO+NO_2)$ cycle through oxidation by ozone ($O_3$). This cycle is key in determining the balance of Ozone ($O_3$) in the atmosphere. The chemical reactions are given by:
\begin{equation}
NO + O_3 \xrightarrow{k_1} NO_2, \quad NO_2 \xrightarrow{k_2} NO + O_3
\label{eq:reactions}
\end{equation}
where the parameters $k_1$ and $k_2$ represent the constant reaction rates throughout the domain. The resulting advection-reaction system of equations can be written as
\begin{eqnarray}
\label{eq:toy_syst} 
\dfrac{\partial \;NO}{\partial t}+u\;\dfrac{\partial \;NO}{\partial x}=-k_1(NO) \;O_3+k_2 \;NO_2 \\
\dfrac{\partial \; NO_2}{\partial t}+u\;\dfrac{\partial \;NO_2}{\partial x}=k_1(NO)\; O_3-k_2 \;NO_2\\
\label{eq:toy_syst2}
\dfrac{\partial \; O_3}{\partial t}+u\;\dfrac{\partial \; O_3}{\partial x}=-k_1(NO) \;O_3+k_2 \;NO_2
\label{eq:toy_syst3}
\end{eqnarray}
where $NO$, $NO_2$, and $O_3$, represent the concentration of each chemical in space and time, and $u$ the constant velocity of the flow (compare with equations (\ref{eq:advec-reac})). \\

The advection and reaction operators commute in this problem (since the advection operator is divergent-free, $\partial u/\partial x=0$, and the chemistry is independent of the location in space), thus, the use of operator splitting approaches should not introduce any error when the exact solutions of the chemistry and advection operators are known. However, when solving numerically the advection operator, with an Eulerian advection scheme, undesired numerical diffusion will cause the numerical advection operator to not commute with the chemistry operator (since nonlinear chemical operators do not commute with diffusion, as shown in \citet{bib:HunVer03}) thus signalling the emergence of operator splitting errors in the numerical solution of equations (\ref{eq:toy_syst})-(\ref{eq:toy_syst3}). \\

This one-dimensional problem is relevant in realistic global 3D simulations since the transport operator is solved utilizing Eulerian numerical schemes, and thus, giving rise to undesired numerical diffusion that will not commute with the time-integration of the chemistry operator. Moreover, in regions in the atmosphere where the flow is near (2D) divergent-free (due to a well stratified atmosphere) and during the night (or day) so that chemistry is independent of space, chemistry and transport operators may commute locally in space and time as in the 1D prototype. \\

In more complicated circumstances, for example in regions of space close to the terminator line (where the day and night boundary is), and in Equatorial regions where convection makes the atmosphere be far from divergent-free conditions, operator splitting errors can be expected to be larger since the advection and chemistry operators will not commute.\\

\subsection{Analytic steady-state solution}

When the chemistry is fast with respect to transport processes, an exact expression can be found for the steady-state solution of system (\ref{eq:toy_syst})-(\ref{eq:toy_syst3}). For example, by choosing $k_1=1000$ and $k_2=2000$, as in \citep{bib:Spo00}, and introducing the non-stiff combined-chemistry operator $\chi = (NO)\; O_3 - 2\; NO_2$, we can represent a stiff (fast) chemistry operator as the quotient $\chi/\epsilon$ for a small parameter $\epsilon$. Equations (\ref{eq:toy_syst}-\ref{eq:toy_syst3}) can be re-written, as suggested in equation (\ref{eq:defn_stiff}), as:
\begin{eqnarray}
\dfrac{\partial \;NO}{\partial t} +u\;\dfrac{\partial \;NO}{\partial x}=-\frac{\chi}{\epsilon}, \label{eq:toy_stiff1} \\ 
\dfrac{\partial \;NO_2}{\partial t}+u\;\dfrac{\partial \;NO_2}{\partial x}=\frac{\chi}{\epsilon}, \\
 \dfrac{\partial \;O_3}{\partial t}+u\;\dfrac{\partial \;O_3}{\partial x}=-\frac{\chi}{\epsilon}.
\label{eq:toy_stiff2}
\end{eqnarray}
Here $\epsilon$ represents the stiffness of the system and is given by the ratio between the slow advection scales and the fast chemistry time scales. For example if $u\sim\mathcal{O}(1)$ and $k_i\sim 10^3$, then $\epsilon\sim10^{-3}$.\\

The expression of the steady-state solution of system is found by introducing the lumped species $NO_x=NO+NO_2$ and $O_x=O_3+NO_2$ (References, sportisse 2000) in order to re-write equations (\ref{eq:toy_stiff1})-(\ref{eq:toy_stiff2}) as:
 \begin{eqnarray}
 \label{eq:toy_stiff_lumped1}
\dfrac{\partial \;NO_x}{\partial t}+u\;\dfrac{\partial \;NO_x}{\partial x}=0, \\
\dfrac{\partial \;O_x}{\partial t}+u\;\dfrac{\partial \;O_x}{\partial x}=0, \\
\dfrac{\partial \;O_3}{\partial t}+u\;\dfrac{\partial \;O_3}{\partial x}=-\frac{\chi}{\epsilon}.
\label{eq:toy_stiff_lumped2}
\end{eqnarray}
In this new form, and denoting  $D/Dt=\partial/\partial t+u \;\partial/\partial x$, it can be seen that the lumped species $NO_x$ and $O_x$ are conserved in time, since 
\[\dfrac{D\;NO_x}{Dt}=0 \quad \text{and} \quad \dfrac{D\;NO_x}{Dt}=0.\] 
As a consequence, for regions where the three species are initially present, the exact asymptotic value of the concentration of all species,  $NO^{\dagger}$, $NO_2^{\dagger}$, and $O_3^{\dagger}$, can be found explicitly as a function of the initial concentration of the lumped species. This is achieved in two steps. First, by expressing the values of the steady state concentrations, $NO^{\dagger}$ and $NO_2^{\dagger}$, as a function of the conserved lumped species as:
\begin{equation}
NO_x(0)=NO^{\dagger}+NO_2^{\dagger}\quad \text{and}\quad O_x(0)=O_3^{\dagger}+NO_2^{\dagger},
\label{eq:lumped_steady}
\end{equation}
and substituting them in equation (\ref{eq:toy_stiff_lumped2}). The system reaches a chemical steady state when $\chi=(NO)\; O_3 - 2\; NO_2=0$, or equivalently when
\begin{equation}
[NO_x(0)- [O_x(0)- O_3^{\dagger}]] O_3^{\dagger} - 2 [O_x(0)- O_3^{\dagger}]=0,
\label{eq:steady_1D}
\end{equation}
which is a second order equation for the steady state of $O_3^{\dagger}$ with solutions given by
\begin{eqnarray}
\label{eq:steady_1D_2} O_3^{\dagger}=&-&\frac{1}{2} \;(2+NO_x(0)-O_x(0))\\ 
&\pm& \frac{1}{2}  \sqrt{(2+NO_x(0)-O_x(0))^2+8 O_x(0)}
\nonumber
\end{eqnarray}
And second, the values of $NO^{\dagger}$, and $NO_2^{\dagger}$ can be found by substituting the (physically relevant) positive solution of (\ref{eq:steady_1D_2}) in equations (\ref{eq:lumped_steady}).  For time scales $\tau$ such that $\tau \gg 1/k$ (for $k=min(k_1, k_2)$), the system will have reached chemical steady-state and from then on, equations (\ref{eq:toy_stiff_lumped1})-(\ref{eq:toy_stiff_lumped2}) (and thus the original system (\ref{eq:toy_syst})-(\ref{eq:toy_syst3})) will behave as a transport-only process propagating the steady-state concentrations with a constant velocity $u$.

\subsection{Numerical experiments}

We chose to solve equations (\ref{eq:toy_syst})-(\ref{eq:toy_syst3}) to simulate the fate of an instantaneous release containing the three chemicals over a $360$ km one-dimensional region. The constant flow velocity was chosen to resemble realistic atmospheric values of $u=10$ m/s. We prescribed a computational spatial domain, $x\in [0,L]$ for $L=3000$ km, so that the plume would stay within the domain for the whole simulation time, $t\in [0,T]$ for $T=10$ hours, and in order to not introduce any errors due to boundary conditions in the numerical advection operator. The values of $k_1=1000$ and $k_2=2000$ were chosen for the stiff chemistry operator. The effective stiffness of the chemistry with respect to the transport is $\mathcal{O}(10^{-2})$ since $u\sim\mathcal{O}(10)$. The initial conditions are given by $NO(x,0)=NO_2(x,0)=O_3(x,0)=p(x)$, where 
\begin{displaymath}
   p(x) = \left\{
     \begin{array}{ll}
       1 & \text{if}\quad x\in[720, 1080]\\
       0 & \text{elsewhere}.
     \end{array}
   \right.
\end{displaymath} 

In a 10-hour simulation time period, the initial release is advected exactly $360$ km to the right, and the concentrations of all species have reached chemical equilibrium. According to expression (\ref{eq:steady_1D_2}), $O_3^{\dagger}=NO^{\dagger}=1.236$, and $NO_2^{\dagger}=0.764$. The exact solution at time $t=T=10$ hours is explicitly given by $O_3(x,T)=NO(x,T)=1.236\times p(x-360)$ and $NO_2(x,T)=0.764\times p(x-360)$. This is our reference solution.\\

\begin{figure}
\centering
\includegraphics[width=.49\textwidth]{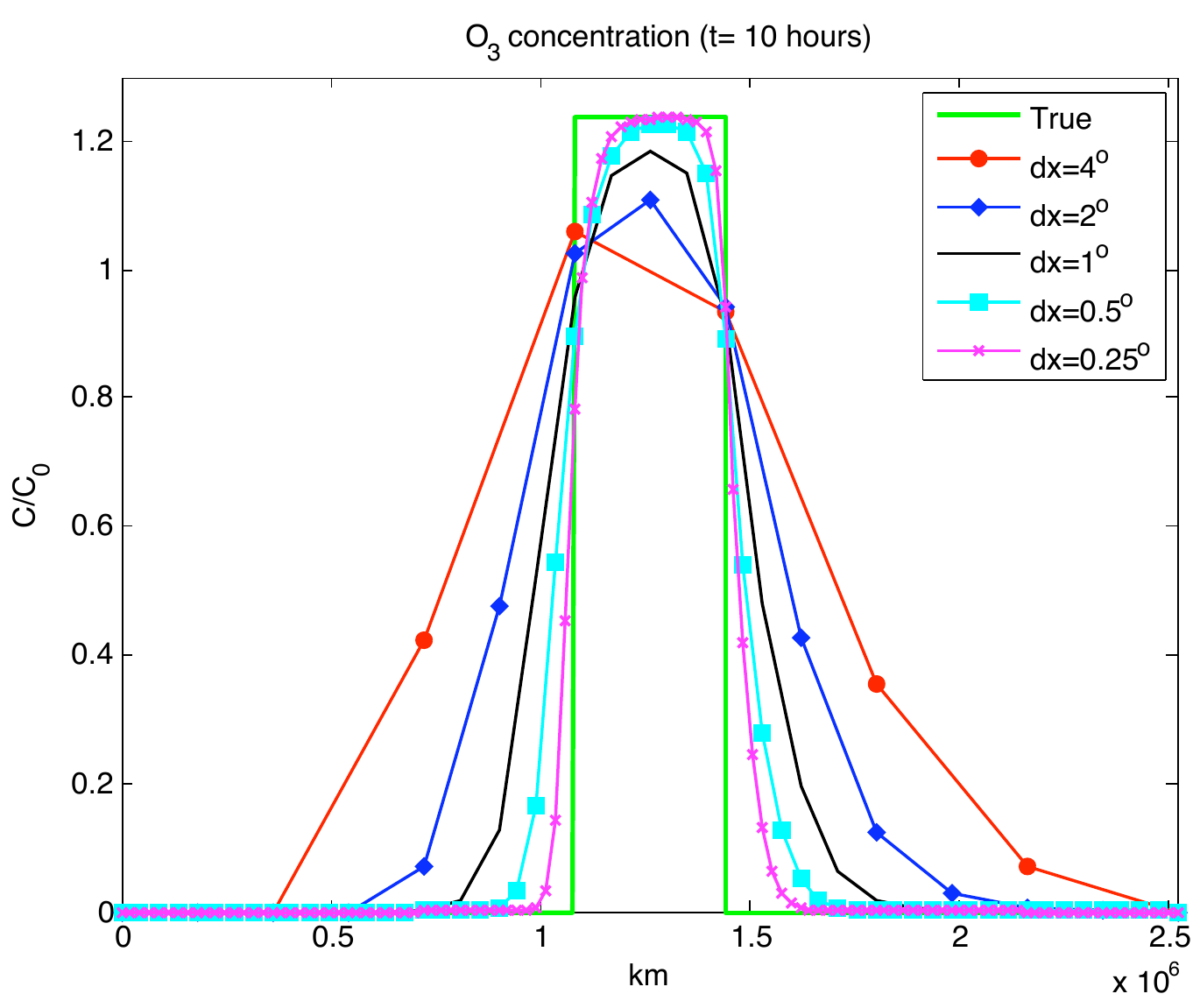}
\includegraphics[width=.49\textwidth]{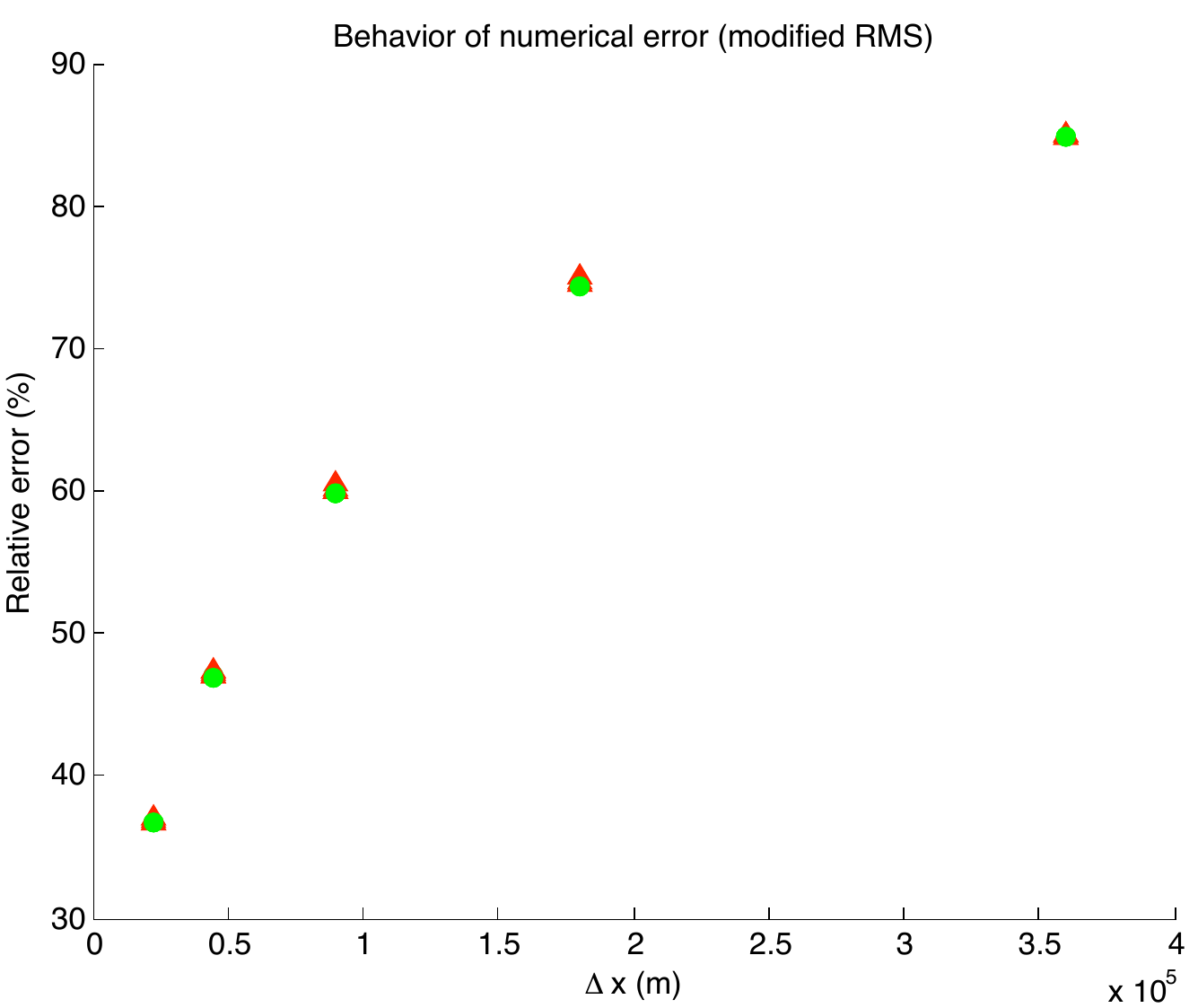}

\includegraphics[width=.49\textwidth]{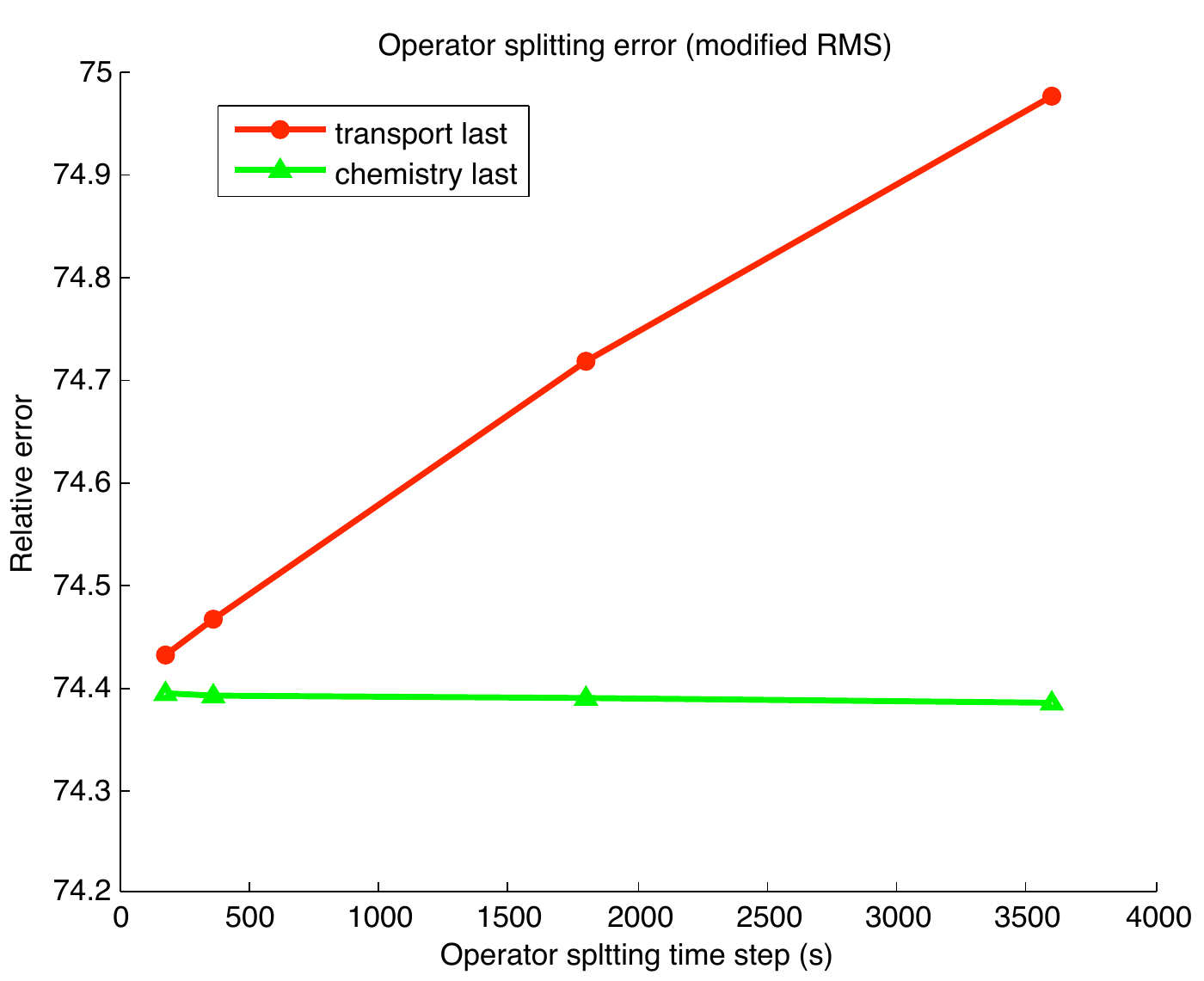}
\caption[Behavior of numerical error ]
{Behavior of numerical error in the one--dimensional transport--chemistry system. The top panel shows the analytical ``true'' and numerical solutions at different grid sizes of the system after a 10-hour simulation time. The middle panel shows the errors relative to the true solution with different grid sizes and operator splitting approaches. The bottom panel shows the behavior of the relative errors (RRMS) from the two operator splitting approaches, for fixed $\Delta x=180$km and different time steps, when compared to the analytic solution.}
\label{fig:ozone_1d}
\end{figure}

For the numerical simulations, we implemented an explicit, second order accurate (in space), one-dimensional advection-scheme based on the Lax-Wendroff method with superbee slope limiters (See \cite{bib:Lev02}, pp 112 for details), and used for the chemistry, the built-in implicit stiff-ODE integrator ode23 from Matlab. In order to minimize contributions to the numerical error, to the first term in inequality (\ref{eq:error_est}), from both the advection scheme and chemistry integrator, we utilized a very small internal advection time step, $\Delta t_{\tau}=90$ seconds, and set the convergence relative-tolerance parameter to $10^{-3}$ in the routine ode23 (it adaptively chooses a small internal time step in order to meet the prescribed $0.1\%$ error convergence criterion). \\

We solved equations (\ref{eq:toy_syst})-(\ref{eq:toy_syst3}) using multiple first order Godunov operator splitting approaches (where transport and chemistry were evaluated in different orders) for multiple operator-splitting time-steps, $\Delta t=180, \;360,\; 1800$ and $3600$ seconds, and for multiple grid sizes $\Delta x=22.5,\;45,\;90, \;180$ and $360$ km (the three largest grid sizes were chosen to resemble spatial resolutions of $4^{\circ}\times 5^{\circ}$, $2^{\circ}\times 2.5^{\circ}$, and $1^{\circ}\times 1.25^{\circ}$, in current 3D global CTMs). The results of these numerical simulations and the exact solution are plotted in the top plot of Figure \ref{fig:ozone_1d}. The numerical solutions corresponding to the multiple operator splitting approaches, for a given value of $\Delta x$, appear as a single curve since their differences were smaller than the line-width chosen for the plot. \\

The quantification of numerical errors was performed using the modified relative root mean square (RRMS), commonly used in 3D atmospheric chemistry simulations given by
\begin{equation}
d_{_{AB}}(C_i)=\sqrt{\frac{1}{M}\displaystyle \sum_{\Omega} \left \vert \dfrac{C_i^A-C_i^B}{C_i^A} \right \vert ^2 }
\label{eq:RMS}
\end{equation}
where $C_i^A$ and $C_i^B$ are the concentrations of species $i$ calculated in simulations $A$ and $B$, respectively, $\Omega$ is the set of grid-boxes where $C_i^A$ exceeds a threshold $a$, 
and $M$ is the number of such grid-boxes. We used $a=$10$^{-4}$, thus neglecting concentrations smaller than $\sim 0.01\%$ with respect to the original concentration. In our one-dimensional experiments, simulation $A$ is the exact solution, and simulation $B$ was one of the multiple Godunov operator splitting approaches. The second plot of Figure \ref{fig:ozone_1d} shows the quantity $d_{_{AB}}=(1/i)\sum_i d_{_{AB}}(C_i)$ for $i=3$ species, for the multiple values of $\Delta t$ and $\Delta x$. In this plot, the red triangles represent simulations where transport was evaluated last, ($\chi-T$), and the green dots where chemistry was evaluated last ($T-\chi$). This plot confirms what is observed in the top plot, {\it i.e.}, the fact that the differences across the multiple operator splitting approaches, for a given $\Delta x$, are very small ($\leq 1\%$).\\

In the bottom plot of Figure \ref{fig:ozone_1d}, we further show the values of the numerical error for the two sequences, $\chi-T$ and $T-\chi$,  for $\Delta x=180$ km, for the multiple values of the operator splitting time-steps. We found this plot to be representative of the behaviour of the numerical error for other values of $\Delta x$. Note that while the differences across the multiple approaches are very small, the interesting mathematical behaviour of the numerical error, discussed in section \ref{sec:op_split}, can be observed. Indeed, the numerical error of the sequences $T-\chi$, where the chemistry (the stiff process) is evaluated last, produce better numerical results than their counter parts $\chi -T$. Moreover, $T-\chi$ sequences appear to be almost insensitive to the magnitude of the operator splitting time-step (the error even seems to grow as $\Delta t\rightarrow 0$ as reported in Sportisse, 2000) making them a preferred choice, since larger operator splitting time steps allow faster computations when exploiting the intrinsic parallelizable nature of the chemistry operator. The quality of results produced by sequences where transport is evaluated last, follows the traditional behaviour of linear analysis where the numerical error decreases as the operator splitting time decreases. Since the magnitude of these first order operator splitting errors was so small, we chose to not implement higher order operator splitting approaches.\\

\begin{figure}
\centering
\includegraphics[width=.49\textwidth]{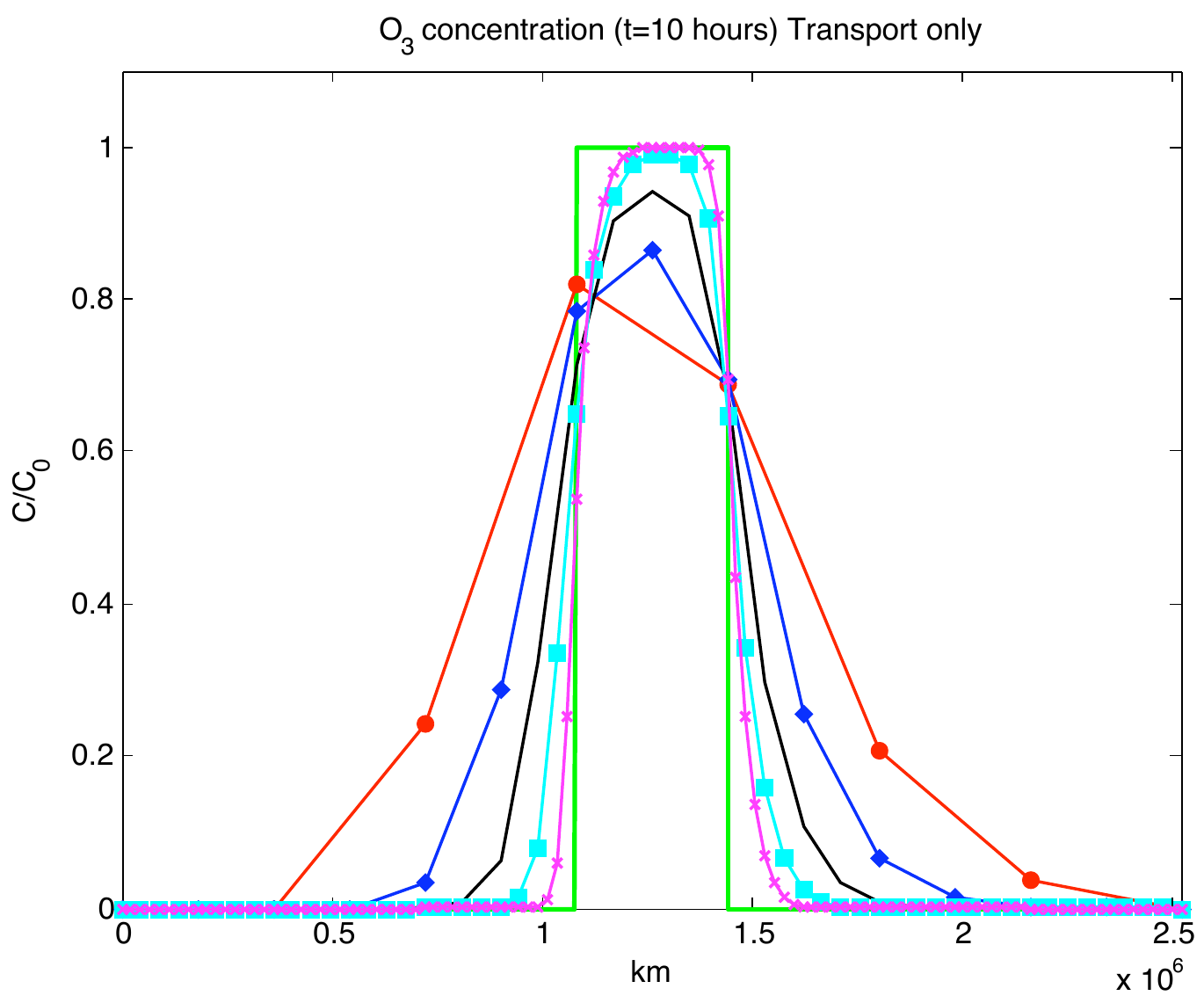}
\includegraphics[width=.49\textwidth]{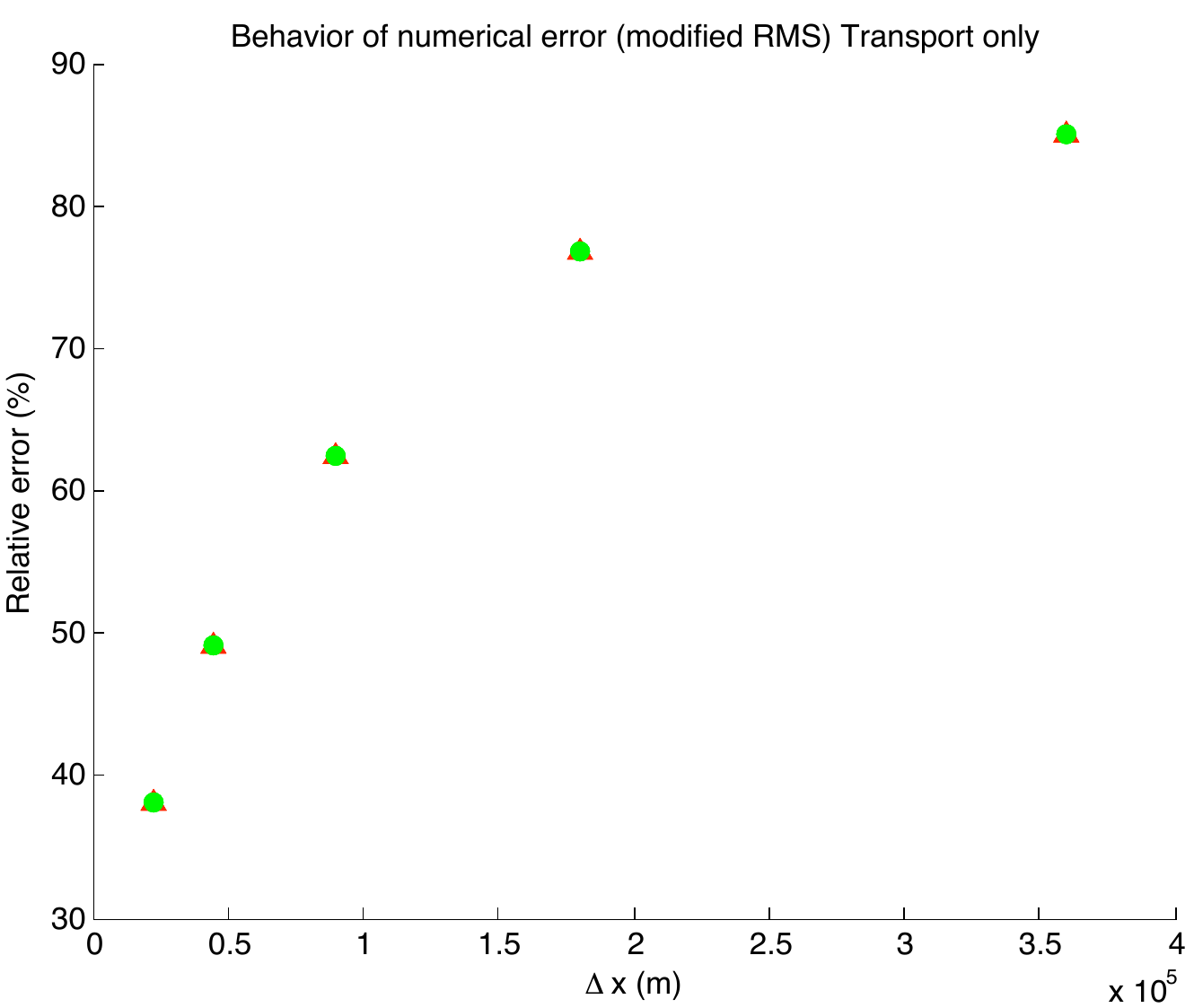}
\caption[Behavior of numerical error, transport only.]
{Behavior of numerical error in the one--dimensional transport--only system. The top panel shows the analytical ``true'' and numerical solutions at different grid sizes of the system after a 10-hour simulation time. The bottom panel shows the errors relative to the true solution with different grid sizes and operator splitting approaches.} 
\label{fig:ozone_1d_transp}
\end{figure}

While the bottom plot of Figure \ref{fig:ozone_1d} shows a clear picture of the magnitude of operator splitting errors ($\leq 1\%$), we performed transport-only simulations in order to verify the magnitude of the numerical errors coming from the numerical advection scheme itself. The results of these simulations are shown in the top plot of Figure \ref{fig:ozone_1d_transp}. Note that while the magnitude of the concentration of $O_3$ in these simulations is exactly one (since no chemistry is present), the numerically simulated profiles, for the different values of $\Delta x$, look very similar to those in the top plot of Figure \ref{fig:ozone_1d}. Indeed, when computing the modified RRMS error associated to these simulations, as shown in the bottom plot of Figure \ref{fig:ozone_1d_transp}, the behaviour of the relative errors resembles the one observed in the middle plot of Figure \ref{fig:ozone_1d}. In short, the numerical errors coming from the choice of operator splitting are eclipsed by the largest component of the numerical error coming from the spatial discretization (second term in inequality (\ref{eq:error_est})) in the numerical advection scheme.\\

Having chosen an initial condition in the shape of a step function in our experiments, caused our second order numerical advection scheme to behave as a first order scheme. Indeed the numerical error decreases close to linearly in our numerical experiments when using the $L^2$-norm instead of the modified RRMS (plot not shown). Estimates for the numerical errors, in the form of an effective numerical diffusion, $D_h$, for 1D first order numerical advection schemes place their value at $D_h\sim u\Delta x$, where $u$ is the mean flow velocity and $\Delta x$ the grid spacing.  In our 1D experiments, these numerical diffusion is of the order $D\sim10^6$ m$^2$/s. Numerical diffusion in 3D global models (\citet{bib:Lin96, bib:San13, bib:Ras07, bib:Wil06, bib:Pis09}) is estimated to be around $10^5-10^6$ m$^2$/s. These 3D estimates place our one-dimensional experiments within a relevant range.

\section{Numerical experiments using GEOS-Chem}
\label{sec:numerical}

Determining the exact magnitude of numerical errors in 3D global CTM simulations in the exact same way we did for our 1D prototype is not possible. This is due to the lack of an analytic expression for the solution to equations (\ref{eq:advec-reac}) in realistic circumstances (time-dependent winds, time-dependent chemistry rates changing throughout the geographic domain due to photolysis, time-dependent emissions).  In order to estimate operator splitting errors in 3D CTMs, we can only compare the output of simulations where everything is kept the same except for the operator splitting sequence and the operator splitting time step. This is the strategy we present in this section, which in combination with the results from our one-dimensional simulations, allowed us to determine upper bounds of operator splitting errors in GEOS-Chem. In order to further understand additional numerical errors, due to the time integration of relevant processes inherently solved with operator splitting approaches and not present in our 1D toy example, we performed multiple additional 3D global simulations. In these simulations, we gradually included inhomogeneous boundary conditions (emission processes) to the time integration (\cite{bib:Spo00, bib:HunVer03}), and vertical processes (convection and dry deposition).\\

GEOS-Chem is a state-of-the-art 3D global Eulerian model of tropospheric chemistry driven by assimilated meteorological observations from the Goddard Earth Observing System (GEOS) of the NASA Global Modeling and Assimilation Ofice (GMAO). The model simulates global tropospheric ozone-NOx-VOC-aerosol chemistry. The full chemical mechanism for the troposphere involves over a hundred species and over three hundred reactions. The  ozone-$NO_x-HO_x-VOC-$aerosol chemical mechanism of GEOS-Chem has been described by \cite{bib:Bey01,bib:Par04} and recently updated by\cite{bib:Mao10}. Details of the chemical reactions and rate constants are reported in the chemical mechanism document (http://acmg.seas.harvard.edu
/geos/wiki\_docs/chemistry
/chemistry \_updates\_v6.pdf). In Figures 4$-$6 the chemical species are arranged in the order of their chemical lifetimes in the atmosphere, from OH ($<$ 1 second) and $NO_x$ ($\sim$1 hour), to CO and $C_2H_6$ (2--3 months). 

The chemical mass balance equations are integrated using a Gear-type solver \citep{bib:Jac95}. Stratospheric chemistry is not explicitly simulated and it instead uses the “Synoz” cross-tropopause ozone flux boundary condition of \citet{bib:McLin00}. The model uses the flux form semi-Langrangian advection scheme of \citet{bib:Lin96}. We used the GEOS-Chem model (v8-02-03) driven by the GEOS-5 data at the 4 x 5 horizontal resolution and 47 levels in the vertical. Detailed descriptions of the model are given by \citep{bib:Bey01} and \citep{bib:Zha11b}. In this study, we initiate the model simulations on January 1, 2005 with model fields from a 6-month spin-up run, and focus on the weekly averaged model results for January 1-7, 2005.


\subsection{Transport and chemistry}

Our strategy consisted of comparing the instantaneous concentration of several chemical species, after multiple one-week long, 4$^\circ$ x 5$^\circ$  horizontal resolution, GEOS-Chem simulations (version v8-02-02), using two versions of the (default) second order Strang operator splitting method given by the sequences: 
\[T(\Delta t/2)  \chi(\Delta t) T (\Delta t/2) \quad \text{and}\quad \chi (\Delta t/2) T(\Delta t) \chi (\Delta t/2)\] 
for different values of the operator-splitting time step $\Delta t$.  These sequences are denoted as $T \chi T$ and $\chi T \chi$ respectively in the subsequent paragraphs. We used $\Delta t= 60, 30, 10, 2$ mins. In all these simulations, transport and chemistry were the only active mechanisms, all other mechanisms were turned off. The inactive mechanisms include: emissions, convective transport, deposition, and planetary boundary layer mixing. Emissions correspond to inhomogenous boundary conditions that are treated numerically as production rates distributed in the boundary layer and solved together in the chemistry operator.\\


We used the modified RRMS (\ref{eq:RMS}) with a threshold $a=$10$^6$ molecules cm$^{-3}$ to quantify the numerical differences in our global simulations. Figure \ref{fig:rms_opsplit_all} shows the relative differences between the reference simulation $\chi T \chi$ with $\Delta t=2$ mins, and the other operator splitting approaches for multiple $\Delta t$'s. Note that the maximum differences across simulations (and species) are of the order of $\sim 10\%$. \\

\begin{figure}
\centering

\includegraphics[width=.49\textwidth]{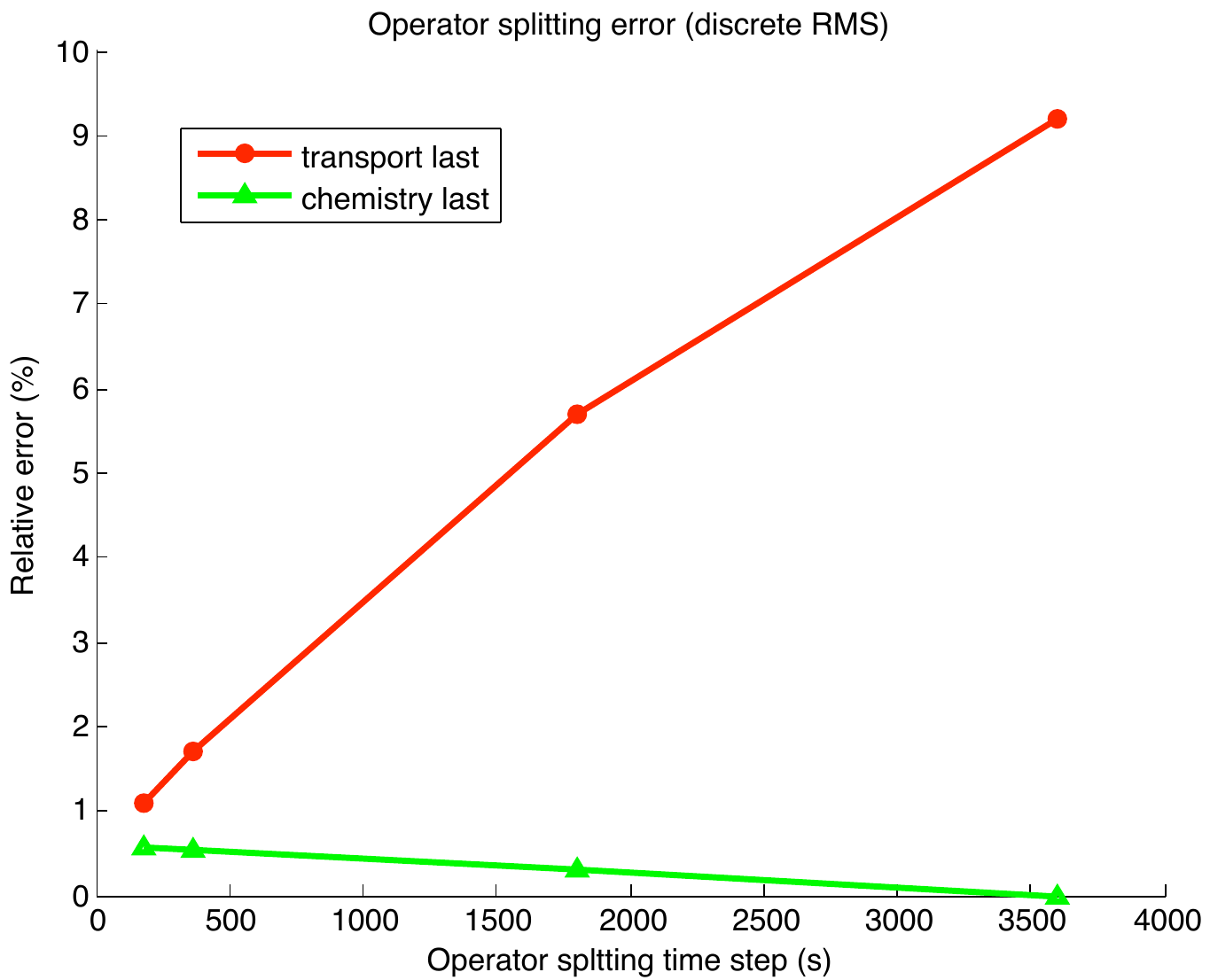}
\caption{Behavior of the relative errors (RRMS) of simulations performed with two different operator splitting approaches ($T-\chi$ and $\chi-T$), fixing $\Delta x=180$ km, for multiples time steps. The reference solution is obtained with the sequence $T-\chi$ for $\Delta t=3600$ seconds.}
\label{fig:1D_opsplit_tch_as_ref}
\end{figure}

Using our one-dimensional prototype and fixing $\Delta x=180$ km, we compared the results of two operator splitting strategies ($T-\chi$ and $\chi-T$) for multiple values of $\Delta t$, with the sequence $T-\chi$ and  $\Delta t=3600$ sec set as a reference. The results are displayed in Figure \ref{fig:1D_opsplit_tch_as_ref}. Note that while the bottom plot of Figure \ref{fig:ozone_1d} shows that operator splitting (relative) errors are less than $ 1\%$ (when comparing to the analytic solution), the relative differences between simulations using alternative operator splitting methods may be as large as $10\%$. This is roughly the same magnitude of the differences observed between the 3D (transport-chemistry) simulations in the top panel of Figure \ref{fig:rms_opsplit_all}.\\

Note that we chose the sequence $\chi T \chi$ with $\Delta t=2$ mins as the reference simulation for our 3D experiments, instead of the sequence $\chi T \chi$ with $\Delta t=60$ that would have been suggested by our 1D experiments (as in Figure \ref{fig:1D_opsplit_tch_as_ref}). The reason for this is shown in the top panel of Figure \ref{fig:rms_op_split_dt}, where we can see that the differences between (transport-chemistry) simulations with different operator splitting sequences but with the same time step, get smaller as $\Delta t$ gets smaller. This behaviour would be expected from a converging operator splitting method where none of the operators is stiff and where the order of evaluation of the operators is not relevant. An alternative explanation could be that the operator splitting errors are very small and what we are observing is the convergence of the time--integration of each operator, separately, as $\Delta t$ gets small. This would suggest that the numerical errors of the time--integration of the transport and the chemistry contribute significantly to the first term (involving $\Delta t$) on the right hand-side of inequality (\ref{eq:error_est}), and should be comparable, in magnitude, to those observed between different operator splitting sequences.\\

In order to investigate this, we plotted the differences between simulations where the only active mechanism was either chemistry or transport, for multiple $\Delta t$'s, while keeping all other parameters exactly the same as in the previous simulations. The results are plotted in Figure 
\ref{fig:rms_chem_only}. These two plots show that indeed the numerical errors arising from the time--integration of each of the operators separately lead to differences of the same magnitude as those observed in the operator splitting simulations. We also observe that the differences get smaller as $\Delta t $ decreases suggesting numerical convergence.  This comparable differences make it hard to disentangle a sharp estimate of the operator splitting in 3D. \\

Note also that in our one dimensional prototype a cleaner analysis was achieved since we chose a smaller internal time step ($\Delta t_{\tau}$= 90 seconds) to integrate the (explicit-in-time) transport operator than the operator splitting time step (180 seconds$\leq\Delta t \leq 60$mins). This choice reduced the contribution to the numerical errors involving $\Delta t$ in inequality (\ref{eq:error_est}) from the transport integration. In order to save computational time in GEOS-Chem (and in most CTMs), however, the time step of the (explicit-in-time) transport scheme is chosen to be equal to the operator splitting time step leading to larger numerical errors. \\

In our one dimensional prototype, the chemistry operator was solved using an adaptive time--integration routine with very tight convergence constraints, thus reducing numerical errors. The time--integration of the chemistry operator in GEOS-Chem uses an adaptive time stepping strategy (\citet{bib:Jac95}) in order to meet convergence requirements (absolute and relative numerical error tolerances) at every user-defined time step. These parameters have been internally set to keep simulation times reasonable while maintaining acceptable numerical accuracy. We kept these settings as they are typically used in global simulations for our numerical experiments. Figure \ref{fig:rms_chem_only} shows the differences between chemistry--only simulations for different user-defined chemistry time-steps. Presumably these errors could be decreased by fine tuning the error tolerances in the time integration routine appropriately, but this approach may increase processing times considerably.\\

Despite all of these numerical issues, we highlight the fact that we can establish an upper limit of about a $10\%$ for the magnitude of operator splitting errors based on the results of our multiple simulations in 3D. Moreover, we show that differences of the single chemical species with largest discrepancies across simulations, Isoprene, are not significant in Figures \ref{fig:isoprene_chem_only}, \ref{fig:isoprene_transp_only}, and \ref{fig:isoprene_op_split_only}, for chemistry only simulations, transport only simulations, and different sequences of operator splitting methods, respectively. From these plots and the results of our one-dimensional prototype, we hypothesize that the operator splitting errors may be much smaller than $10\%$. \\

We also highlight the fact that we did not pursue further efforts to show that the sequences evaluating the chemistry at the end of the time step in 3D compare better with observations, since our one-dimensional prototype, as well as multiple studies in global CTMs \citep{bib:Ras07, bib:Pra08, bib:San13}, suggest that the numerical errors associated to the transport integration, at current spatial resolutions, are significantly larger than those observed in operator splitting methods. In addition, uncertainties  in emission fields and deposition mechanisms may pose further difficulties in addressing this question. In our one-dimensional proto-type, subsequent reductions in the spatial resolution lead to significant improvements in the accuracy of the numerical solution globally (for any operator splitting sequence). Whereas a better choice of operator splitting (where chemistry is evaluated last) leads to a very modest improvement at a given spatial resolution $\Delta x$.

\subsection{Boundary conditions and vertical processes}

Other important processes in 3D simulations are integrated in time using operator splitting strategies. As noted in \cite{bib:Spo00} and \cite{bib:HunVer03} the time integration of inhomogeneous boundary conditions, such as emission processes in global simulations, using operator splitting strategies may lead to considerable numerical errors. Additionally, the time integration of vertical processes such as convection and deposition using operator splitting may also lead to important numerical errors. 

In order to investigate the magnitude of numerical errors due to these processes, we performed additional 3D simulations that gradually included inhomogeneous boundary conditions (emissions) and vertical process. In other words, aside from the 3D ``transport-chemistry'' simulations discussed in the previous sections, we performed simulations with (i) ``transport, chemistry, and emissions'' and simulations with (ii) ``transport, chemistry, emissions, convective transport and deposition''. When emissions are included, they are integrated within the chemistry solver, using the chemistry time step. Convective transport and deposition are solved using the standard setting of GEOS-Chem, which integrate these two processes (sequentially) during the chemistry time step. 

\begin{figure}
\centering

\includegraphics[width=.49\textwidth]{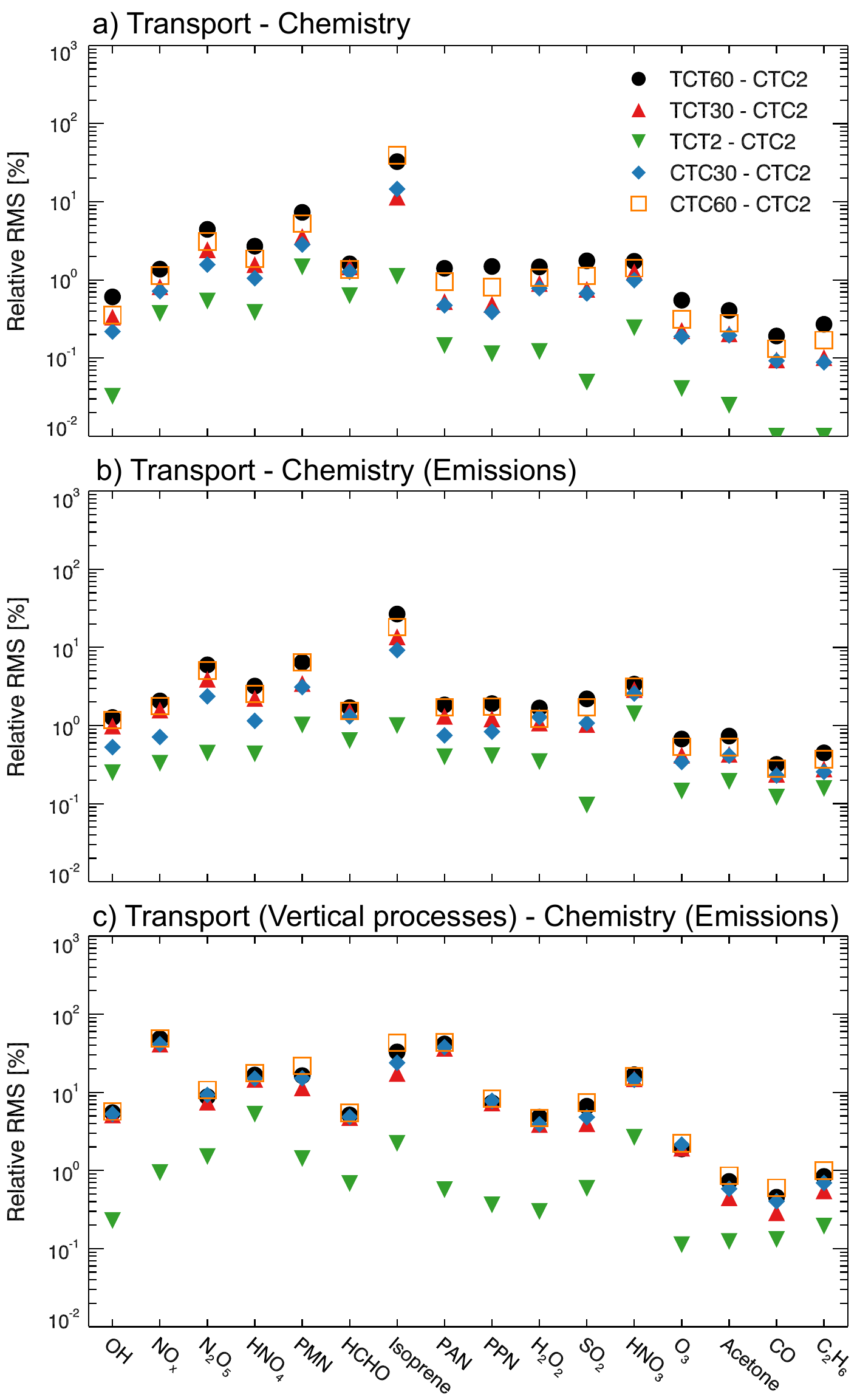}

\caption{Behavior of numerical error in the GEOS-Chem 3-D model simulations. Here TCT denotes Transport-Chemistry-Transport, CTC denotes Chemistry-Transport-Chemistry, and the numbers denotes operator splitting time steps in minutes. Relative RMS relative to the CTC2 model simulation are shown for different chemical species with lifetimes ranging from seconds $(OH)$ to months $(CO, C_2H_6)$. Active processes in these simulations are as follows: Transport and chemistry (top panel); Transport, chemistry and emissions (middle panel); Trasnport, chemistry, emissions, convective transport and deposition (bottom panel). } 
\label{fig:rms_opsplit_all}
\end{figure}

\begin{figure}
\centering

\includegraphics[width=.49\textwidth]{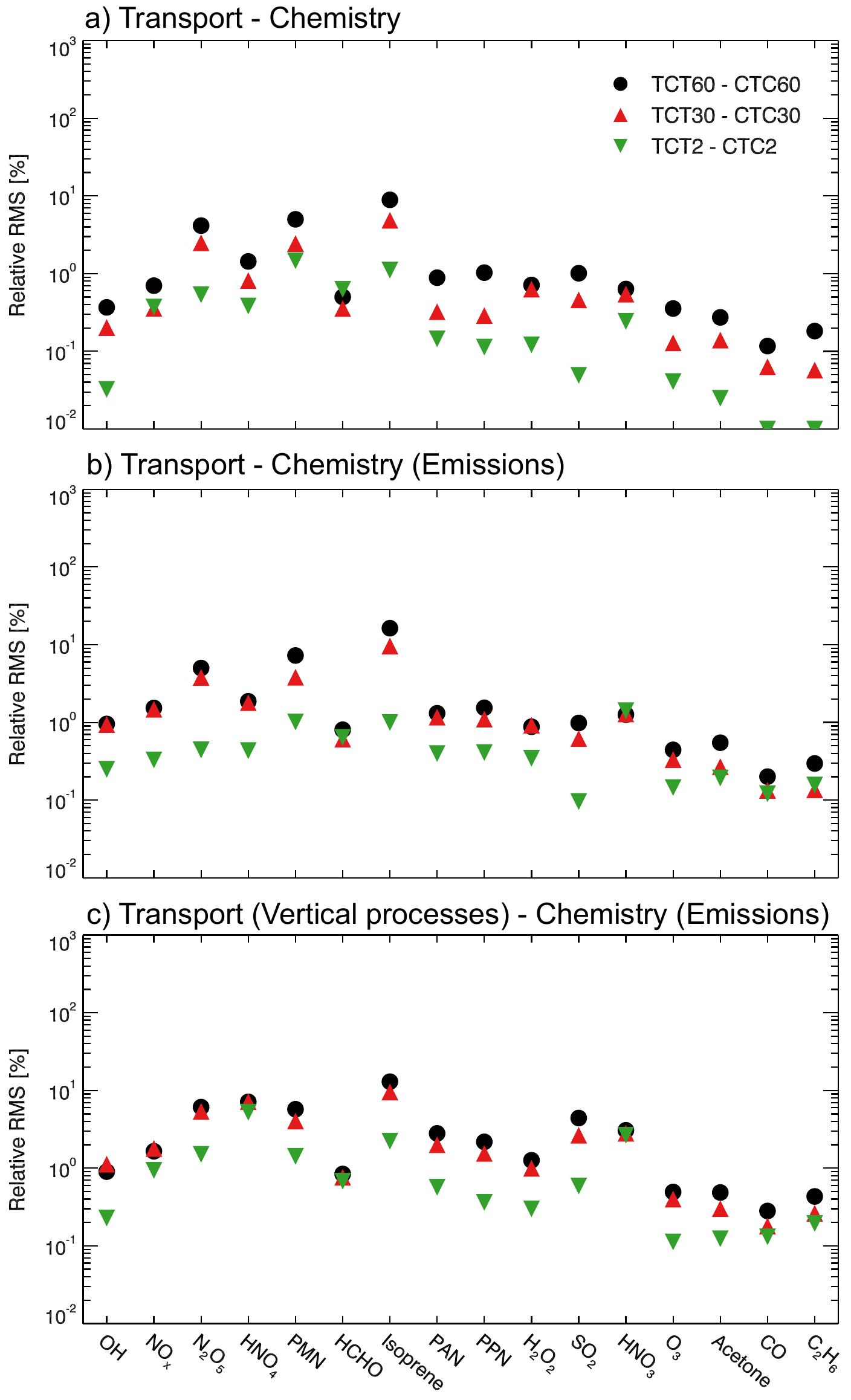}

\caption{Behavior of numerical error in the GEOS-Chem 3-D model simulations. Here TCT denotes Transport-Chemistry-Transport, CTC denotes Chemistry-Transport-Chemistry, and the numbers denotes operator splitting time steps in minutes. Relative RMS for different operator splitting approaches for fixed time steps: $\Delta t=2,30, 60$ mins.Active processes in these simulations are as follows: Transport and chemistry (top panel); Transport, chemistry and emissions (middle panel); Trasnport, chemistry, emissions, convective transport and deposition (bottom panel).}
\label{fig:rms_op_split_dt}
\end{figure}

\begin{figure}
\centering

\includegraphics[width=.49\textwidth]{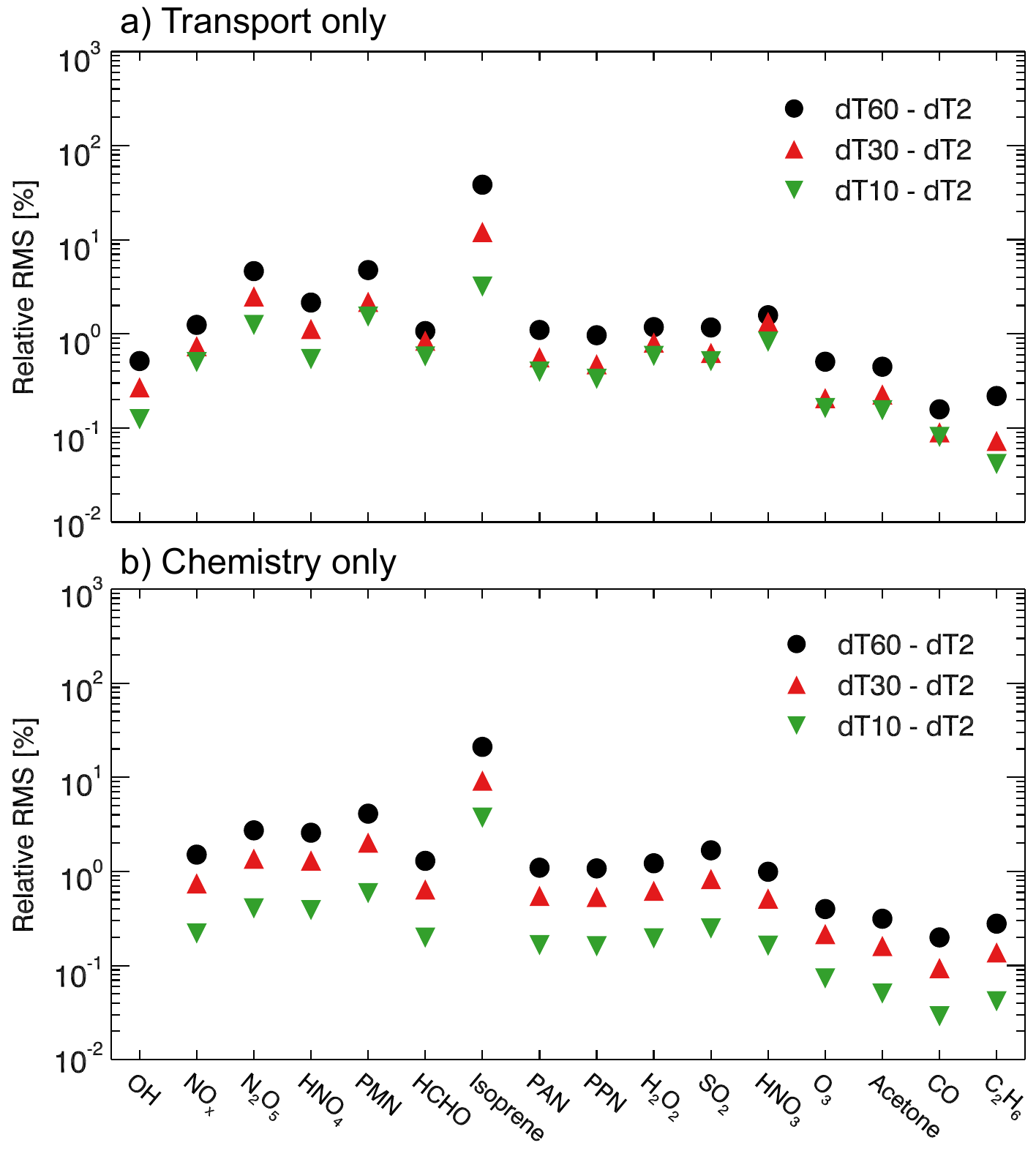}

\caption{Behavior of numerical error in the GEOS-Chem 3-D model simulations. Relative RMS for transport--only (top panel) and chemistry only (bottom) simulations using different time steps: $\Delta t=2,30, 60$ mins.}
\label{fig:rms_chem_only}
\end{figure}

\begin{figure}
\centering

\includegraphics[width=.49\textwidth]{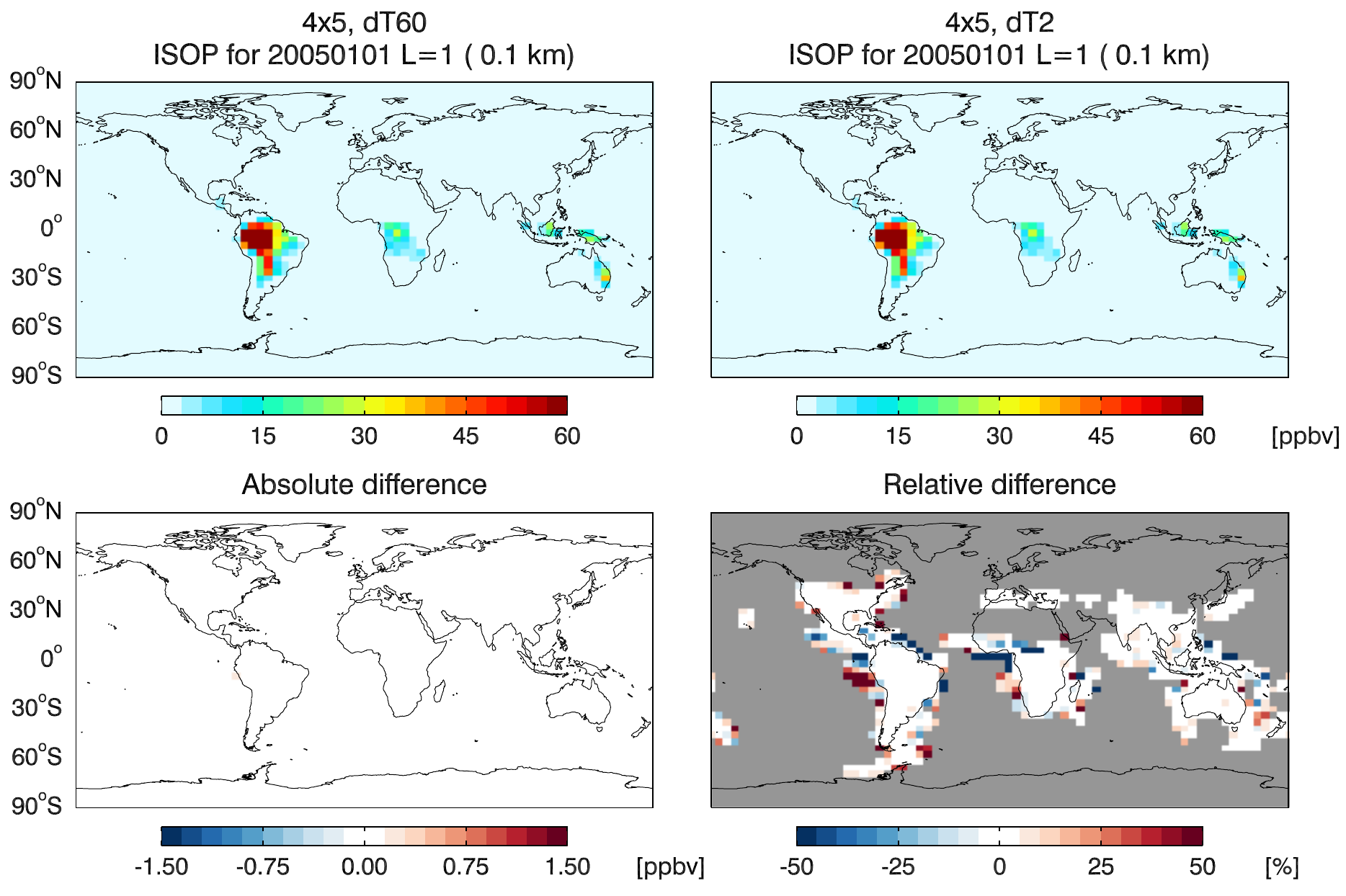}

\caption{Comparison of isoprene concentrations using different time steps for GEOS-Chem transport--only simulations. Isoprene concentrations at the surface level from the model simulation with time step of 60 minutes (top-left panel) are compared to the model simulation with time step of 2 minutes (top-right panel). Absolute (bottom-left) and relative differences (bottom-right) are also shown.}
\label{fig:isoprene_transp_only}
\end{figure}

\begin{figure}
\centering

\includegraphics[width=.49\textwidth]{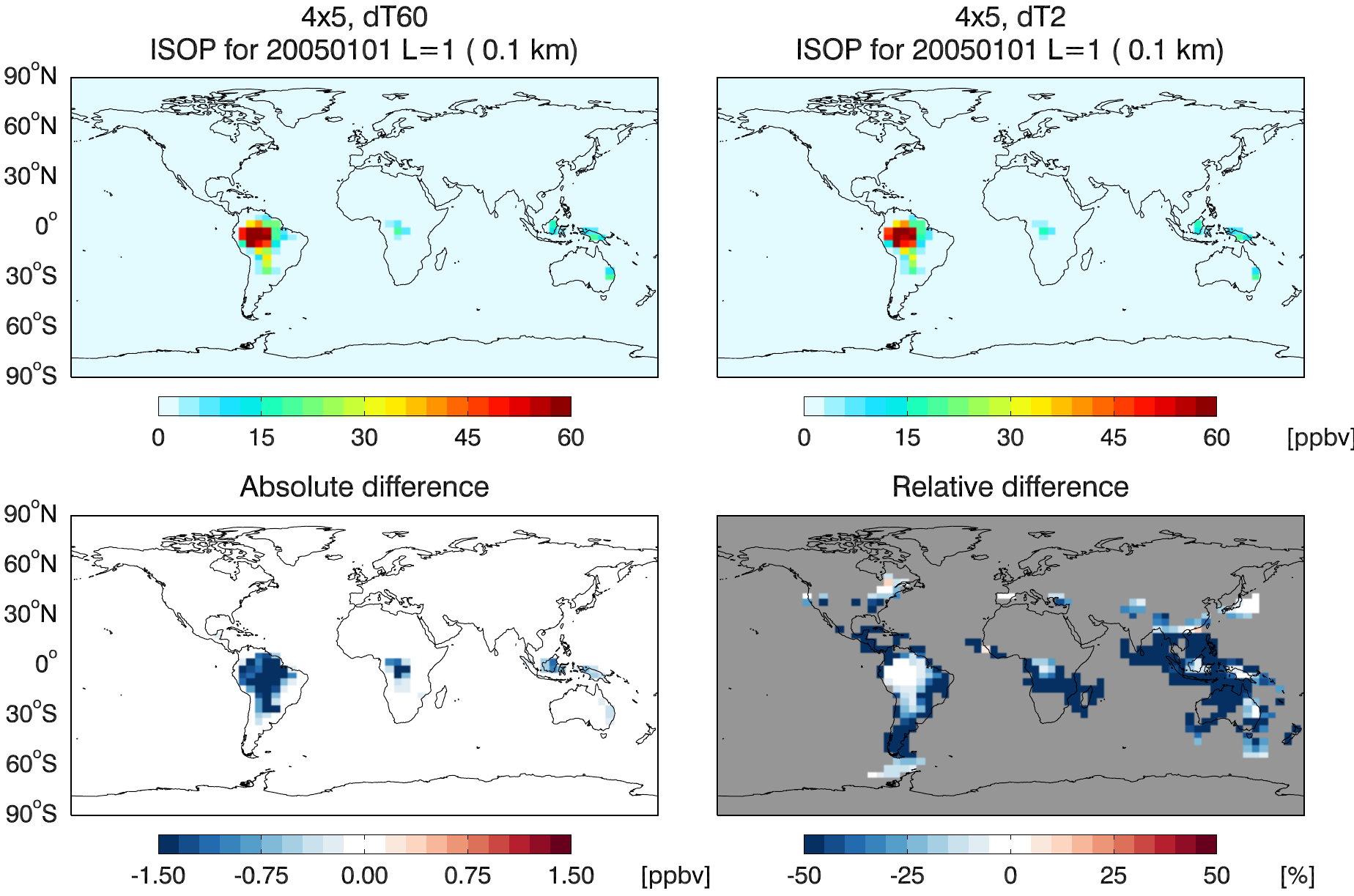}

\caption{Comparison of isoprene concentrations using different time steps for GEOS-Chem chemistry--only simulations. Isoprene concentrations at the surface level from the model simulation with time step of 60 minutes (top-left panel) are compared to the model simulation with time step of 2 minutes (top-right panel). Absolute (bottom-left) and relative differences (bottom-right) are also shown.}
\label{fig:isoprene_chem_only}
\end{figure}

\begin{figure}
\centering

\includegraphics[width=.49\textwidth]{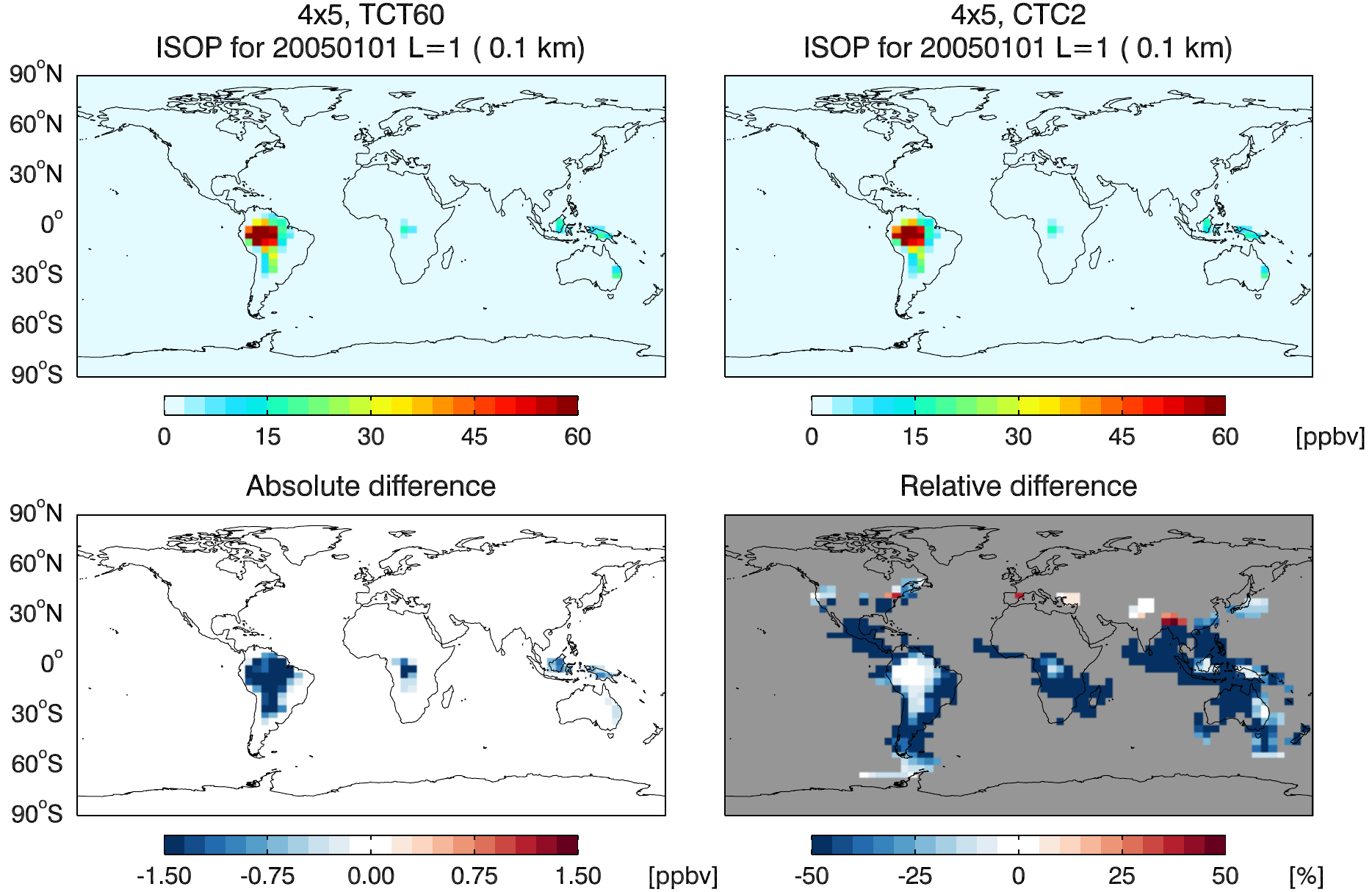}

\caption{Comparison of isoprene concentrations using Transport-Chemistry-Transport (time step of 60 minutes) versus Chemistry-Transport-Chemistry (time step of 2 minutes).}
\label{fig:isoprene_op_split_only}
\end{figure}

The differences between these two sets of simulations, using the same methodology explained in the previous section, are plotted in the two lower panels of Figures \ref{fig:rms_opsplit_all} and \ref{fig:rms_op_split_dt}. As these Figures show, the additional numerical errors coming from the inclusion of inhomogenous boundary conditions (emissions) are significant. Indeed, the differences between the simulations that include ``transport, chemistry, and emissions'' are roughly double the magnitude of the differences between the simulations that include only ``transport-chemistry'' for different operator splitting strategies. The incorporation of convective transport and deposition to the simulations does increase the differences between simulation, mainly when changes in time steps are large, as shown in the bottom panel of Figure \ref{fig:rms_opsplit_all}. When time-steps are fixed and operator splitting approaches are different, these vertical processes do not seem to lead to larger differences in the different simulations.

\section{Conclusions and Future work}
We have presented a way to characterize operator splitting errors in the context of atmospheric chemistry modeling. Our approach numerically extends one--dimensional linear results to non-linear 1D and 3D cases. These numerical findings are relevant to global atmospheric chemistry modeling. Our findings suggest that stiff operators should be evaluated lastly in operator splitting methodologies. This results is consistent with the linear results presented in \citep{bib:Spo00}, and previous studies in numerical weather prediction \citep{bib:Dub05} . Differences of approximately $10\%$ across species are found when comparing the outputs of global simulations using different operator splitting approaches, using multiple splitting time steps. This, in combination with our one-dimensional results, suggests that operator splitting errors do not exceed $10\%$ relative errors in global simulations. We show also, that in current spatial resolutions, the numerical diffusion errors introduced in global atmospheric chemistry models eclipse errors emerging from operator splitting techniques.

\subsection{Future work}
Future studies should identify whether operator splitting strategies that evaluate fast dynamics operators lastly in global simulations lead to simulations that improve the match between simulations and observations. Further exploration is also required regarding the effect of different operator splitting strategies in the time integration of of the governing equations of aerosol dynamics and different choices of boundary layer mixing schemes. Additional ``toy-tests'' that should be explored in order to further understand numerical errors introduced by different operator splitting strategies include those discussed in \citep{bib:Lau14, bib:Pud06}. Finally, nuances between operator splitting approaches in Eulerian and Semi-Lagrangian transport schemes should be more deeply investigated \citep{bib:Pud97}.

%
%

\section*{Acknowledgements}
MS and LZ would like to thank the technical assistance provided by Claire Carouge. MS would like to thank Jonathan Pines for his involvement in the exploratory phases of this project. This work was partially funded by the National Natural Science Foundation of China (41205103)


\begin{thebibliography}{00}

\bibitem[Allen et al.(1996)]{bib:All96}
Allen, Dale J., Richard B. Rood, Anne M. Thompson, and Robert D. Hudson. Three-dimensional radon 222 calculations using assimilated meteorological data and a convective mixing algorithm. Journal of Geophysical Research: Atmospheres 101 (1996), no. D3: pp. 6871--6881.doi: 10.1029/95JD03408

\bibitem[Bey et al.(2001)]{bib:Bey01}
Bey I., D. J. Jacob, R. M. Yantosca, J. A. Logan, B. Field, A. M. Fiore, Q. Li, H. Liu, L. J. Mickley, and M. Schultz. Global modeling of tropospheric chemistry with assimilated meteorology: Model description and evaluation, Journal of Geophysical Research, 106 (2001), D19: pp. 23073--23095. doi:10.1029/2001jd000807.

\bibitem[Brenner and Scott(2008)]{bib:Bre08}
Brenner, S. C., and Scott R. The mathematical theory of finite element methods. Vol. 15. Springer Science and Business Media, (2008). doi: 10.1007/978-0-387-75934-0

\bibitem[Dubal et al.(2005)]{bib:Dub05}
Dubal, Mark, Nigel Wood, and Andrew Staniforth. Mixed parallel-sequential-split schemes for time-stepping multiple physical parameterizations. Monthly weather review 133 (2005), no. 4: pp. 989--1002. doi:10.1175/MWR2893.1

\bibitem[Enting(2002)]{bib:Ent02} 
Enting, I. G. Inverse Problems in Atmospheric Constituent Transport. Cambridge University Press, (2002). doi:10.1017/CBO9780511535741 

\bibitem[Ern and Guermond(2004)]{bib:Ern04}
Ern, A., and Guermond, J.L. Theory and practice of finite elements, vol. 159 of Applied Mathematical Sciences, (2004). doi:10.1007/978-1-4757-4355-5

\bibitem[Guo and Babuska (1986)]{bib:Guo86}
Guo, B., and I. Babuska. The hp version of the finite element method. Computational Mechanics 1.1 (1986): pp. 21--41. doi: 10.1007/BF00298636

%
%
%


\bibitem[Hundsdorfer and Verwer(2003)]{bib:HunVer03}
Hundsdorfer W. and Verwer J.G. Numerical Solution of Time--Dependent 
Advection--Diffusion--Reaction Equations, Springer Series in Computational Mathematics (2003), 33, Springer. doi: 10.1007/978-3-662-09017-6

\bibitem[Iserles(2009)]{bib:Ise96}
A. Iserles.: A first course in the numerical analysis of differential equations, Cambridge University Press (2009), 44. doi: 10.1007/978-81-322-1835-7


\bibitem[Jacobson(1995)]{bib:Jac95}
Jacobson M. Z. Computation of Global Photochemistry with SMVGEAR-II, Atmospheric Environment, 29(18) (1995), pp. 2541--2546. doi:10.1016/1352-2310(95)00194-4


\bibitem[Lansen and Verwer(1999)]{bib:LanVer99}
Lanser D., Verwer J.G. Analysis of operator splitting for advection–-diffusion–-reaction problems from air pollution modelling, Journal of Computational and Applied Mathematics. 111 (1999), pp. 201–-216. doi: 10.1016/S0377-0427(99)00143-0

\bibitem[Lauritzen(2014)]{bib:Lau14}
Lauritzen, P. H., Conley, A. J., Lamarque, J. F., Vitt, F., \& Taylor, M. A., The terminator toy-chemistry test: a simple tool to assess errors in transport schemes, Geoscientific Model Development Discussions, 7 (6)(2014), pp. 8769--8804. doi:10.5194/gmd-8-1299-2015

\bibitem[LeVeque(2002)]{bib:Lev02}
LeVeque, R. J. Finite volume methods for hyperbolic problems. Vol. 31. Cambridge university press (2002). ISBN-13: 978-0521009249

\bibitem[Lin and Rood(1996)]{bib:Lin96}
Lin, S.J, and R. B. Rood. Multidimensional flux-form semi-Lagrangian transport schemes. Monthly Weather Review 124(9)(1996), pp. 2046--2070. doi 10.1175/1520-0493(1996)124

\bibitem[Lowe and Tomlin(2000)]{bib:LowTom00}
Lowe, R. and Tomlin, A. Low-dimensional manifolds and reduced chemical models for tropospheric chemistry simulations, Atmospheric Environment, 34 (2000), pp. 2425--2436. doi:10.1016/S1352-2310(99)00447-1

\bibitem[Mao et al(2010)]{bib:Mao10}
J. Mao, D. J. Jacob, M. J. Evans, J. R. Olson, X. Ren, W. H. Brune, J. M. St. Clair, J. D. Crounse, K. M. Spencer, M. R. Beaver, P. O. Wennberg, M. J. Cubison, J. L. Jimenez, A. Fried, P. Weibring, J. G. Walega, S. R. Hall, A. J. Weinheimer, R. C. Cohen, G. Chen, J. H. Crawford,C. McNaughton,A. D. Clarke,L. Jaeglé, J. A. Fisher, R. M. Yantosca, P. Le Sager, and C. Carouge. Chemistry of hydrogen oxide radicals (HOx) in the Arctic troposphere in spring, Atmospheric Chemistry and Physics, 10 (2010), pp. 5823--5838. doi:10.5194/acp-10-5823-2010

\bibitem[McLinden et al.(2000)]{bib:McLin00}
McLinden, C. A., Olsen, S. C., Hannegan, B., Wild, O., Prather, M. J., \& Sundet, J. Stratospheric ozone in 3-D models: A simple chemistry and the cross-tropopause flux, Journal of Geophysical Research, 105(D11) (200), pp. 14653-–14666. doi: 10.1029/2000JD900124

\bibitem[Park et al(2004)]{bib:Par04}
Park, R. J., D. J. Jacob, B. D. Field, R. M. Yantosca, and M. Chin. Natural and transboundary pollution influences on sulfate--nitrate--ammonium aerosols in the United States: implications for policy, Journal of Geophysical Research, 109 (2004), D15204, doi: 10.1029/2003JD004473

\bibitem[Pisso et al.(2009)]{bib:Pis09}
Pisso, I., Elsa Real, Kathy S. Law, B. Legras, N. Bousserez, J. L. Attié, and H. Schlager. Estimation of mixing in the troposphere from lagrangian trace gas reconstructions during long--range pollution plume transport. Journal of Geophysical Research: Atmospheres, 114 (D19) (2009). doi:10.1029/2008JD011289

\bibitem[Prather et al.(2008)]{bib:Pra08}
Prather, M. J., et al. Quantifying errors in trace species transport modeling. Proceedings of the National Academy of Sciences 105(50) (2008). pp. 19617--19621. doi: 10.1073/pnas.0806541106


\bibitem[Pudykiewicz et al.(1997)]{bib:Pud97}
Pudykiewicz J., A. Kallaur, and P. Smolarkiewicz, Semi-Lagrangian Modeling of Tropospheric Ozone, Tellus 49B (1997), pp. 231–-248. doi: 10.1034/j.1600-0889.49.issue3.1.x

\bibitem[Pudykiewicz(2006)]{bib:Pud06}
Pudykiewicz J., Numerical Solution of the Reaction--Advection--Diffusion Equation on the Sphere, Journal of Computational Physics, 213(1) (2006), pp. 358-–390. doi:10.1016/j.jcp.2005.08.021


\bibitem[Rastigeyev et al.(2007)]{bib:Ras07}
Rastigejev, Y., M. P. Brenner, and D. J. Jacob. Spatial reduction algorithm for atmospheric chemical transport models. Proceedings of the National Academy of Sciences 104(35) (2007), pp 13875--13880. doi: 10.1073/pnas.0705649104


\bibitem[Santillana et al.(2010)]{bib:San10}
M. Santillana, P. Le Sager, D. J. Jacob, and M. P. Brenner. An adaptive reduction algorithm for efficient chemical calculations in global atmospheric chemistry models, Atmospheric Environment 44 (35) (2010), pp. 4426--4431.         doi:10.1016/j.atmosenv.2010.07.044

\bibitem[Santillana(2013)]{bib:San13}
M. Santillana. Quantifying the loss of information in source attribution problems using the adjoint method in global models of atmospheric chemical transport. arXiv.org: 1311.6315, (2013)

\bibitem[Sportisse(2000)]{bib:Spo00}
B. Sportisse. An Analysis of Operator Splitting Techniques in the Stiff Case, Journal of Computational Physics, 161 (2000), pp. 140–-168. doi: 10.1006/jcph.2000.6495

\bibitem[Sportisse(2007)]{bib:Spo07}
B. Sportisse. A review of current issues in air pollution modeling and simulation. Computational Geosciences 11.2 (2007) pp. 159--181. doi: 10.1007/s10596-006-9036-4

\bibitem[Strang(1968)]{bib:Str68}
G. Strang. On the construction and comparison of difference schemes. SIAM Journal of Numerical Analysis 5.3 (1968), pp. 506–-517. doi: 10.1137/0705041

\bibitem[Wild and Prather(2006)]{bib:Wil06}
Wild, O., and M. J. Prather. Global tropospheric ozone modeling: Quantifying errors due to grid resolution. Journal of Geophysical Research: Atmospheres 111.D11305 (2006). doi:10.1029/2005JD006605

\bibitem[Zhang et al(2011)]{bib:Zha11}
H. Zhang, J. C. Linford, A. Sandu, and R. Sander. Chemical Mechanism Solvers in Air Quality Models, Atmosphere 2(3) (2011), pp. 510--532. doi:10.3390/atmos2030510

\bibitem[Zhang et al(2011)b]{bib:Zha11b}
L. Zhang, D. J. Jacob, N. V. Downey, D. A. Wood, D. Blewitt, C. C. Carouge, A. van Donkelaar, D. BA Jones, L. T. Murray, and Y. Wang. Improved estimate of the policy-relevant background ozone in the United States using the GEOS-Chem global model with 1/2 × 2/3 horizontal resolution over North America, Atmospheric Environment 45 (2011), pp. 6769–-6776. doi:10.1016/j.atmosenv.2011.07.054


\end{thebibliography}
\end{document}